%% file: Main.tex
\documentclass[12pt,a4paper]{amsart}
\usepackage{graphicx} 
\usepackage[utf8]{inputenc} 
\usepackage{latexsym,amssymb,amsmath,amsfonts,amsthm,mathtools,mathrsfs}
\usepackage{amscd}
\numberwithin{equation}{section}
\renewcommand{\raggedright}{\justifying}							
\usepackage{helvet}																		
\usepackage{calligra}
\usepackage{enumitem}
\usepackage{url}
\usepackage{filecontents}
\input{comandos}

\begin{document}

\title{Hyperbolic manifolds with a large number of systoles}


\author{Cayo D\'oria}

\author{Emanoel M.~S. Freire} 
\thanks{The second author was supported by a CAPES research grant.}

\author{Plinio G. P. Murillo}

\address{
Cayo Dória \newline
Universidade Federal de Sergipe, Departamento de Matemática. \newline
Av. Marcelo D\'eda Chagas s/n, 49100-000. São Cristóvão - SE, Brazil}
\email{cayo@mat.ufs.br}
\address{
	Emanoel M. S. Freire \newline
	IMPA\\
	Estrada Dona Castorina, 110\\
	22460-320 Rio de Janeiro, Brazil}
\email{emanoel.m.s.freire@gmail.com}
\address{
Plinio G. P. Murillo\newline
Universidade Federal Fluminense, Instituto de Matemática e Estatística.\\
Rua Prof. Marcos Waldemar de Freitas Reis, s/n, Bloco H, Campus do Gragoatã, 24210-201. Niterói-RJ, Brazil.}
\email{pliniom@id.uff.br}

\begin{abstract} 
In this article, for any $n\geq 4$ we construct a sequence of compact hyperbolic $n$-manifolds $\{M_i\}$ with number of systoles at least as $\vol(M_i)^{1+\frac{1}{3n(n+1)}-\epsilon}$ for any $\epsilon>0$. In dimension 3, the bound is improved to $\vol(M_i)^{\frac{4}{3}-\epsilon}$. These results generalize previous work of Schmutz for $n=2$, and Dória-Murillo for $n=3$ to higher dimensions. 
\end{abstract}

\maketitle

\section{Introduction}
\input{Sec1_intro}
\section{Preliminaries}\label{preliminar}
\input{Sec2_preli}
\section{Length inequality in $\Spin(1,n)$}
\input{Sec3_length_inequality}
\section{Arithmetic subgroups of $\Spin(1,n)$}\label{congrugroup}
\input{Sec4_5_Congruence_subgroups}
\section{Proof of Theorem \ref{firstresult} and Theorem \ref{maintheorem}}\label{proofofmainresult}
\input{Sec6_Proofs}
\section{The three-dimensional case}\label{sec2}
\input{Sec7_dimension_three}
\bibliographystyle{amsalpha}
\bibliography{ref}

\end{document}

%% file: comandos.tex
\newtheorem{theorem}{Theorem}[section]
\newtheorem{lemma}[theorem]{Lemma}
\newtheorem{corollary}[theorem]{Corollary}

\newtheorem{theoremletra}{Theorem}

\newtheorem{prop}[theorem]{Proposition}

\theoremstyle{definition}

\newtheorem{question}{Question}

\newcommand{\bl}{\begin{lemma}}
\newcommand{\el}{\end{lemma}}
\newcommand{\bprop}{\begin{prop}}
\newcommand{\eprop}{\end{prop}}
\newcommand{\bt}{\begin{theorem}}
\newcommand{\et}{\end{theorem}}
\newcommand{\bcor}{\begin{corollary}}
\newcommand{\ecor}{\end{corollary}}
\newcommand{\Spin}{\mathrm{Spin}}

\theoremstyle{remark}
\newtheorem{remark}{Remark}
\newtheorem*{notation*}{Notation}

\usepackage[color=red!40]{todonotes}

\newcommand{\inv}{^{-1}}

\newcommand{\Hyp}{\mathbb{H}^3}

\newcommand{\N}{\mathrm{N}}
\newcommand{\Z}{\mathbb{Z}}
\newcommand{\Q}{\mathbb{Q}}
\newcommand{\R}{\mathbb{R}}
\newcommand{\HS}{\mathbb{H}}
\newcommand{\CC}{\mathbb{C}}
\newcommand{\SL}{\mathrm{SL}}
\newcommand{\PSL}{\mathrm{PSL}}
\newcommand{\PSLC}{\mathrm{PSL}(2,\mathbb{C})}
\newcommand{\sys}{\mathrm{sys}}
\newcommand{\tr}{\mathrm{tr}}

\newcommand{\K}{\mathrm{Kiss}}

\newcommand{\area}{\mathrm{area}}
\newcommand{\vol}{\mathrm{vol}}
\newcommand{\Isom}{\mathrm{Isom}}
\newcommand{\hol}{\mathrm{hol}}
\newcommand{\Arg}{\mathrm{Arg}}
\newcommand{\ok}{\mathcal{O}_k}
\newcommand{\nm}{\mathcal{N}}

\newcommand{\Kiss}{\mathrm{Kiss}}

\usepackage{hyperref}
\DeclareFontFamily{T1}{calligra}{}
\DeclareFontShape{T1}{calligra}{m}{n}{<->s*[1.2]callig15}{}

\usepackage{enumitem}
\usepackage{mathtools}

%% file: Sec1_intro.tex
A \emph{hyperbolic $n$-manifold} is a $n$-dimensional manifold without boundary, equipped with a complete Riemannian metric of constant curvature $-1$. Any hyperbolic manifold is isometric to a quotient space $M= \Gamma \backslash \mathbb{H}^n$, where $\mathbb{H}^n$ is the hyperbolic $n$-space and $\Gamma$ is a torsion-free discrete group of isometries of $\mathbb{H}^n$. In recent years a lot of progress has been made in the study of hyperbolic manifolds with extremal properties. For example, such spaces with minimal volume \cite{Bel14}, minimal diameter \cite{BCP21}, large systolic length \cite{Mur19} and large kissing number (see below the definitions of systolic length and kissing number).
 In this article we are interested in hyperbolic manifolds with large  kissing number, and their relation with its systolic length and volume.

It is well known that any finite volume hyperbolic manifold $M$ contains closed geodesics, and there exists at least one of minimal length, which is called a \emph{systole} of $M$. The length of any systole of $M$ will be called the \emph{systolic length} of \(M\) and denoted by $\sys(M)$ .

We define the \emph{kissing number} $\K(M)$ of $M$ as the number of systoles of $M$. This is a well-defined invariant since, in negative curvature there are finitely many closed geodesics with the same length. This terminology was introduced by Schmutz Schaller, and it was inspired by the classical kissing number of lattices arising in sphere packings (see \cite{Sch96a},\cite{Sch96b},\cite{Sch95}).

In general, it follows from a classical result of Anosov \cite{Anosov83} that a generic Riemannian manifold has at most one systole. For a closed hyperbolic $n$-manifold $M$, it is possible to bound $\Kiss(M)$ by above in terms of $\sys(M)$ or $\vol(M)$. 
It follows from works by Keen \cite{Kee74} and Buser \cite{Bus80} that there exist constants \(C_n,D_n>0\) depending only on \(n\) such that if \(\sys(M) \leq C_n\), then
\begin{equation}\label{C.27.06.22}
\Kiss(M) \leq D_n \vol(M). 
\end{equation}

Recently, Bourque and Petri  proved an upper bound of $\Kiss(M)$ which holds for any systolic length (see \cite[Theorem 1]{BP19}). More precisely,
\begin{equation}\label{kissystole}
\K(M) \leq A_n \vol(M) \frac{\exp\left( \frac{n-1}{2}\sys(M)\right)}{\sys(M)},
\end{equation}
for some constant \(A_n>0\) which depends only on \(n.\) Since $\frac{\sinh(x)}{x}\rightarrow 1$ whenever $x\rightarrow0$, we can unify the inequalities \eqref{C.27.06.22} and \eqref{kissystole} in 
 \begin{equation}\label{kissinguppergeneral}
     \K(M) \leq A'_n \vol(M) \frac{\sinh \left(\frac{n-1}{2}\sys(M) \right)}{\sys(M)}
 \end{equation}
 for some constant \(A'_n>0\) which depends only on \(n\) \footnote{We thank the referee for providing us this remark.}.
However, if $\sys(M)$ is large, Inequality \eqref{kissystole} implies that

\begin{equation}\label{upperbounds}
\K(M) \leq B_n\frac{\vol(M)^2}{\log(1+\vol(M))},
\end{equation}
(see \cite[Corollary 1.2]{BP19}). In dimension 2, this result was previously established in $2013$ by Parlier (see \cite{Parlier13}). In \cite{FP16}, similar upper bounds were established for non-compact hyperbolic surfaces of finite area. Whether a version of \eqref{kissystole} and \eqref{upperbounds} holds for non-compact finite volume hyperbolic manifolds in $n\geq 3$ remains open.

These restrictions for $\Kiss(M)$, and the aforementioned result by Anosov motivated us to study the following question formulated in \cite{Petri20}: Let $n \geq 2$ and
	\[K_n(v)=\max\{\Kiss(M) \mid M \mbox{ is a hyperbolic $n$-manifold of } \vol(M)\leq v \} \]
 \begin{question} \label{questione}
 How does $K_n(v)$ grows as a function of $v$?
\end{question}
Although this question is independent of the size of $\sys(M)$, it is interesting to understand $\Kiss(M)$ according to whether $\sys(M)$ is small or large. The first result in this article shows that there exist a sequence of closed hyperbolic manifolds $N_i$ with bounded systolic length, and $\Kiss(N_i)$ growing at least as linearly with $\vol(N_i)$. 
\begin{theoremletra}\label{firstresult}
	Let $n \geq 2$ and $A>0$. Then, there exist positive real numbers $B,C$ with $A\leq B$, and a sequence $\{N_i\}$ of closed hyperbolic $n$-manifolds such that 

\begin{itemize}
\item $A\leq \sys(N_i) \leq B$ for all $i>0$,

\item $\vol(N_i) \to \infty$

\item $\K(N_i) \geq C \cdot \vol(N_i)$.
\end{itemize}

	In particular,
	\[\limsup\limits_{v \to \infty} \frac{\log K_n(v)}{\log v} \geq 1 \]
	for any $n \geq 2$.
\end{theoremletra}
This result shows that the exponent one of $\vol(M)$ in \eqref{kissystole} is the best possible that we can obtain. The manifolds in Theorem \ref{firstresult} are obtained by taking cyclic covers of a fix hyperbolic manifold with positive first Betti number. The use of this technique is well known in spectral geometry  (see e.g \cite{Randol74}). The proof of Theorem \eqref{firstresult} is given in Section \ref{proofofmainresult}.

The main part of the article is devoted to give an answer to Question \ref{questione} independently on the size of the systoles.  For $n=2$, it follows from results by Schmutz in \cite{Sch95} that
\begin{equation}\label{eq:K_2 Schmutz}
\limsup\limits_{v \to \infty} \dfrac{\log K_2(v)}{\log v} \geq 1+\frac{1}{3}.
\end{equation}

To prove this result, in [op. cit.] the author constructed a sequence $S_i$ of closed (also non-compact of finite area) hyperbolic surfaces with large kissing number as congruence coverings of a fixed arithmetic hyperbolic surface. It is worth to note that the surfaces \({S_i}\) also satisfy $$\sys(S_i) \sim \frac{4}{3}\log (\area(S_i))\xrightarrow{i\rightarrow\infty}\infty.$$ 


More generally, if a sequence $M_i$ of non-diffeomorphic closed hyperbolic $n$-manifolds has $\Kiss(M_i)$ growing super linearly in $\vol(M_i)$ (i.e. $\Kiss(M_i) \geq C \vol(M_i)^{1+\varepsilon}$ for some constants $C, \varepsilon>0$), then $\sys(M_i)$ grows logarithmic in $\vol(M_i)$. Indeed, it follows from \eqref{kissinguppergeneral} since $\vol(M_i) \to \infty$, and $\dfrac{\sinh x}{x} >M$ implies $x>\log(M)$ for any $M>0$. Hence 
\begin{equation}\label{eq:sysvolgeneral}
  \sys(M_i) \geq \frac{2\varepsilon}{n-1} \log(\vol(M_i))+\frac{2}{n-1}\log\left(\frac{2C}{A_{n}(n-1)}\right).
\end{equation}

In \cite{Mur19}, the third author showed that congruence coverings of closed arithmetic hyperbolic $n$-manifold of the first type have systole with length growing logarithmically in the volume, and determined the precise growth ratio. It is then natural to investigate the kissing number of such manifolds, and to ask whether they can provide a version of \eqref{eq:K_2 Schmutz} in higher dimension. In this direction we obtain the following

\begin{theoremletra}\label{maintheorem}
For any $n \geq 2$, there exists a  compact arithmetic hyperbolic n-manifold of the first type $M$, and a sequence of congruence coverings $M_j \to M$ of arbitrarily large degree such that
\begin{equation} \label{mainpartofmaintheorem}
\K(M_j) \geq  C\ \frac{\vol(M_j)^{1+\frac{1}{3n(n+1)}}}{\log(\vol(M_j))} 
\end{equation}
for some constant $C>0$ independent of $M_j$.
In particular,
	\[\limsup\limits_{v \to \infty} \frac{\log K_n(v)}{\log v} \geq 1 + \frac{1}{3n(n+1)}\]
for any $n \geq 2.$
\end{theoremletra}

Although the systolic length of congruence coverings of arithmetic manifolds is large, exhibiting the systoles in these spaces is a much more delicate problem. We overcome this by constructing $M_j$ containing a totally geodesic surface $S_j$ whose systoles are also systoles of $M_j$.

\begin{theoremletra}\label{th:embeded_surface}
	For any $n \geq 3$, the manifold $M$ obtained in Theorem \ref{maintheorem}  contains a closed totally geodesic surface $S$ such that for any $j$, the congruence coverings $M_j \to M$ contains a congruence covering $S_j \to S$ satisfying  $\sys(S_j)=\sys(M_j).$ 
\end{theoremletra}

This is the key part of the proof of Theorem \ref{maintheorem}  (see Corollary \ref{equalsystoles}).  The result then follows from an argument of high multiplicity inspired by \cite{Sch95}. The proof is presented in Section \ref{proofofmainresult}.

The last part of the article reserves a special attention to dimension $3$. In \cite{DM21}, the first and last authors constructed non-compact hyperbolic manifolds $N_i$ satisfying  \[\frac{\log \K(N_i)}{\log \vol(N_i)} \gtrsim 1+\frac{4}{27}.\] 

In this case, the manifolds $N_i$ are congruence coverings of Bianchi orbifolds. Once more, the technique includes ideas by Schmutz, but a new length-trace relation was needed, and also a result on averages of class numbers of imaginary binary forms proved by Sarnak. These results are no longer available in the compact setting. Instead, we have now been able to use the holonomy of closed geodesics. We are able to construct closed arithmetic hyperbolic $3$-manifolds with a large number of systoles using the relation between length and trace of $2 \times 2$ matrices, and a result on the equidistribution of closed geodesics with holonomy in prescribed intervals proved by Sarnak and Wakayama in \cite{SW99} (see also \cite{Marg14}). We finish the article with the proof of the following theorem in Section \ref{sec2}.

\begin{theoremletra}\label{secondresult}
 There exists a sequence $\{M_j\}$ of compact arithmetic hyperbolic $3$-manifold with $\vol(M_j)$  going to infinity such that
 \[\K( M_j) \geq C \frac{\vol(M_j)^{4/3}}{\log(\vol(M_j))}\]
where $C>0$ is a universal constant.

 \end{theoremletra}

\textbf{Comments and an open question.} We recall that two positive sequences \((a_j)\) and \((b_j)\) satisfy the relation \(a_j \gtrsim b_j\) (resp. \(a_j\lesssim b_j\))  when for any \(\delta>0\) there exists \(j_0\) such that \(\frac{a_j}{b_j} \geq (1-\delta)\) (resp. \(\frac{a_j}{b_j}\leq 1-\delta\)) for \(j>j_0\). Hence, the sequences satisfy \(a_j \sim b_j\) if and only if \(a_j \lesssim b_j\) and \(a_j \gtrsim b_j. \) A natural question arising from Theorem \ref{maintheorem} is the following: 
\begin{question}\label{questiontwo}
Is there a universal $\varepsilon>0$ such that for any $n\geq 2$, there is a sequence of closed hyperbolic $n$-manifolds $M_j$ with  \(\vol(M_j) \to \infty \) and \[\K(M_j) \gtrsim \frac{\vol(M_j)^{1+ \varepsilon}}{\log(\vol(M_j))}  ~?\]
\end{question}

We have already noticed that this would imply 
\begin{equation}\label{eq:sys_question}
\sys(M_j) \gtrsim \frac{2 \varepsilon}{n-1}\log(\vol(M_j))
\end{equation} (see Inequality \eqref{eq:sysvolgeneral}). From the Appendix of \cite{Mur19}, for any $n\geq 2$ there is a sequence of compact arithmetic hyperbolic $n$-manifolds $M_i$ with $\vol(M_i)\rightarrow\infty$ such that
\begin{equation}\label{C091221}
\sys(M_i) \sim \frac{8}{n(n+1)}\log(\vol(M_i)),
\end{equation}
and the hypothetical bound \eqref{eq:sys_question} would be considerably larger than the growth in \eqref{C091221}. While writing this article the authors do not know any improvement for \eqref{C091221}.

\vskip 10pt

\noindent{\textbf{Acknowledgements.}}
We would like to thank the referees for their comments and suggestions which substantially improved the previous version of this article.

%% file: Sec2_preli.tex
\subsection{Hyperbolic Manifolds}

The hyperbolic $n$-space is the complete simply connected $n$-dimensional Riemannian manifold with constant sectional curvature equal to $-1$. The hyperboloid model of the hyperbolic $n$-space is given by
$$\mathbb{H}^n=\{x \in \mathbb{R}^{n+1};\  -x_0^2+x_1^2+\cdots + x_n^2=-1, \ x_0>0\}$$
with the metric $ds^2=-dx_0^2+dx_1^2+\cdots+dx_{n}^2$.

The identity component $\mathrm{SO}(n,1)^{\circ}$ of the Lie group $\mathrm{SO}(n,1)$ 
is isomorphic to the orientation preserving isometries of $\mathbb{H}^{n}$. Given a lattice $\Gamma\subset\mathrm{Isom}^{+}(\mathbb{H}^n)$, i.e, a discrete subgroup having finite covolume with respect to the Haar measure, the associated quotient space $M=\Gamma\backslash\mathbb{H}^n$ is a finite volume  {\it hyperbolic orbifold}. It is a manifold whenever $\Gamma$ is torsion-free. 

\subsection{Systole and Kissing number}
We recall that an element $\gamma$ in $\mathrm{Isom}^{+}(\mathbb{H}^n)$ is called
\begin{itemize}
\item {\it elliptic}, if it has a fixed point on $\mathbb{H}^n$.
\item {\it parabolic}, if it has exactly one fixed point in $\partial \mathbb{H}^n$.
\item {\it loxodromic}, if it has  two fixed points in $\partial \mathbb{H}^n$.
\end{itemize}

The displacement at $x \in \mathbb{H}^n$ of a loxodromic element $\gamma$ is defined by $l(\gamma,x):=d(x,\gamma x)$. The {\it displacement} of $\gamma$ (also called {\it translation length}) is defined by
$$l(\gamma):= \mathrm{inf}_{x\in \mathbb{H}^n} l(\gamma,x).$$

A \textit{systole} of a hyperbolic orbifold $M=\Gamma\backslash\mathbb{H}^{n}$ is any closed geodesic of shortest length in \(M\), and its length is denoted
by $\sys(M)$. The {\it kissing number} $\mathrm{Kiss}(M)$ is defined as the number of free homotopy classes of oriented closed geodesics in $M$ that realize $\sys(M)$.
It is well known that when \(M\) is a manifold, there exists a one-to-one correspondence between parametrized closed geodesic in $M$ up to unit speed reparametrization and conjugacy classes of loxodromic elements in $\Gamma$. Moreover, the length of the closed geodesic corresponds to the translation length of the loxodromic element. This relation allows to study $\sys(M)$ and $\mathrm{Kiss}(M)$ through the number of conjugacy classes of loxodromic elements in $\Gamma$ that realize $\sys(M)$ as their translation length.

\subsection{Clifford algebras and the Spin group}

In this section we will recall the construction of spinor groups. These are the algebraic groups associated to fundamental groups of  arithmetic hyperbolic orbifolds. For further details, we refer the reader to \cite[Section 2]{Mur19}, and the references therein.
 
Let $k$ be a field with char~$k\neq 2$, $E$ a $n$-dimensional vector space over $k$, and $f$ a non-degenerate quadratic form on $E$ with associated bilinear form $\Phi$. The \textit{Clifford algebra of $f$}, denoted by $\mathscr{C}(f,k)$, is a unitary associative algebra over $k$ given by the quotient $\mathscr{C}(f,k)=T(E)/\mathfrak{a}_{f}$, where $T(E)$ denotes the tensor algebra of $E$, and $\mathfrak{a}_{f}$ is the two-sided ideal of $T(E)$ generated by the elements $x\otimes y +y\otimes x - 2\Phi(x,y)$. If we choose an orthogonal basis $e_1,\dots,e_{n}$ of $E$ with respect to $f,$ we have in $\mathscr{C}(f,k)$ the relations $e_{v}^2=f(e_{v})$ and $e_{v}e_{\mu}=-e_{\mu}e_{v}$ for $\mu,v\in\{1,\ldots,n\}, \ \mu \neq v.$ 
Let $\mathscr{P}_{n}$ be the power set of $\{1,...,n\}$. For $M=\{\mu_1,\dots,\mu_\nu\}\in \mathscr{P}_{n}$ with $\mu_1<\dots<\mu_\nu,$ we define $e_M=e_{\mu_1}\cdot...\cdot e_{\mu_\nu}$, where we adopt the convention $e_{\emptyset}=1$. The $2^n$ elements $e_M$, $M\in\mathscr{P}_n$ determine a basis for $\mathscr{C}(f,k)$. Hence, any \(s \in \mathscr{C}(f,k)\) is written uniquely as $$s=s_\R \cdot 1 +\sum_{\mathclap{M \in \mathscr{P}_{n}, ~ M \neq \emptyset}}s_{M}e_{M},$$ with \(s_{\R}, s_M \in k\). We call the coefficient \(s_\R\) as the \textit{real part} of $s$. We identify $k$ with $k\cdot 1$, and $E$ with the $k$-linear subspace in $\mathscr{C}(f,k)$ generated by $e_1, e_2,\ldots, e_n$. 
The algebra $\mathscr{C}(f,k)$ has an anti-involution $*$. On the elements $e_M$ this map acts by $e_{M}^{*}=(-1)^{\nu(\nu-1)/2}\ e_{M}$, where \(\nu=|M|\). 
The span of the elements $e_M$ with $|M|$ even is a subalgebra $\mathscr{C}^{+}(f)$ of $\mathscr{C}(f),$ called the \textit{even Clifford subalgebra of $f$}. The spin group of $f$ is defined as
\begin{equation*}
\Spin_{f}(k):=\left\{ s \in  \mathscr{C} ^{+}(f,k)\ \Big| \ sEs^*\subseteq E, \ ss^*=1 \right\}.    
\end{equation*}

In the case $k=\mathbb{R}$, $E=\mathbb{R}^{n+1}$ and $f=-x_0^2+x_1^2+\cdots+x_n^2$, the corresponding spin group is denoted by $\Spin(n,1)$. For an element $s\in\Spin(n,1)$ the linear map $\varphi_{s}:E\to E$ given by $\varphi_{s}(x)=sxs^{*}$ lies in $\mathrm{SO}(n,1)^{\circ}$, and the map $s\to\varphi_{s}$ is a two-sheeted covering of $\mathrm{SO}(n,1)^{\circ}$, with kernel $\{\pm 1 \}$.
Since the image of a lattice under a finite covering map is also a lattice, in order to produce hyperbolic orbifolds we will implicitly contruct lattices in $\Spin(n,1)$ and project them to $\mathrm{SO}(n,1)^{\circ}$. Furthermore, we will abuse terminology saying that an element \(s \in \Spin(n,1)\) is elliptic (resp. parabolic or loxodromic) if $\varphi_{s}$ is elliptic (resp. parabolic or loxodromic).

%% file: Sec3_length_inequality.tex
In order to obtain explicitly the systole of a closed manifold $M=\Gamma \backslash \mathbb{H}^n$, we need to find a hyperbolic element $\gamma_0\in \Gamma$ such that $\ell(\gamma_0) \leq \ell(\gamma)$ for any nontrivial element $\gamma \in \Gamma$. As it has been observed in \cite{Mur19}, it is useful to estimate the displacement of hyperbolic elements using information about the real part of elements in the spin group. Now we will give a more precise version of this relation. The main tool will be the connection between spin and Vahlen groups established by Elstrodt, Grunewald and Mennicke \cite{EGM87}, and a characterization of the translation length by Waterman \cite{Wat}.

\begin{prop}\label{prop:realpart-length}
For any loxodromic element $r \in \Spin(1,n)$ we have that $$\ell(r)\geq 2\cosh^{-1}(|r_{\mathbb{R}}|).$$
\end{prop}

\begin{proof}

Let us consider $k=\mathbb{R}$ and the quadratic forms

\begin{align*}
Q_0(y_0,y_1,y_2) &=y_0^2-y_1^2-y_2^2\\
Q(x_1,x_2,\ldots,x_{n-2}) & = -x_1^2-\ldots-x_{n-2}^2\\
\tilde{Q}(y_0,y_1,y_2,x_1,\ldots,x_{n-2}) &=Q_0(y_0,y_1,y_2) + Q(x_1,\ldots,x_{n-2}).
\end{align*}

It is clear that $\Spin_{\tilde{Q}}(\mathbb{R})=\Spin(1,n)$. By \cite[Theorem 4.1]{EGM87}, any element $r\in\Spin(1,n)$ is the image under a group isomorphism $\psi$ of a $2\times 2$ matrix with coefficients in the Clifford algebra $\mathscr{C}(Q)$. More precisely, and following the notation in [\textit{loc. cit}] a direct computation shows that

\begin{equation}\label{eq:dif_cliff_algebras}
\begin{aligned}
r&=\psi\left(\begin{pmatrix}
a & b\\
c & d \\
\end{pmatrix} \right) \\
&=\dot{a}\frac{1}{2}(1+f_0f_1)+\dot{b}\frac{1}{2}(f_0f_2-f_1f_2)+\dot{c}\frac{1}{2}(f_0f_2+f_1f_2)+\dot{d}\frac{1}{2}(1-f_0f_1)\\
\end{aligned} 
\end{equation}
where $a, b, c, d \in \mathscr{C}(Q)$. 
Here $f_0,$ $f_1$ and $f_2$ denote a basis of ortho\-gonal elements for $Q_0$, and $\cdot$ denotes the $K$-algebra homomorphism from $\mathscr{C}(Q)$ to $\mathscr{C}(\tilde{Q})$ given by

$$\dot{x}=\left(\sum_{M\in \mathscr{P}_{n}}x_{M}e_{M}\right)^{\cdot}=\sum_{M\in \mathscr{P}_{n}}x_{M}(f_0f_1f_2)^{\xi_{M}}e_{M}, $$

where
\begin{equation*}
  \xi_{M}=\left\{ \begin{array}{rll}
0 \ \text{for} \ |M|\equiv 0\ \text{mod}\ 2  \\
\\
1 \ \text{for} \ |M|\equiv 1\ \text{mod}\ 2, \\
\end{array} \right.   
\end{equation*}
(see \cite[Section 2]{EGM87}). The important point for us is that $f_0, f_1, f_2$ are not elements in the center of $\mathscr{C}(\tilde{Q})$, and the real part of $\dot{x}$ is equal to that of $x$, for any $x\in\mathscr{C}(Q)$. Looking at \eqref{eq:dif_cliff_algebras}, this implies that the real part of $r$ is determined by the real part of $a+d$. More precisely 
\begin{equation}\label{eq:length_parte_real_1}
r_\mathbb{R}= (a+d)_{\mathbb{R}}/2.
\end{equation}

On the other hand, by \cite[Lemma 14]{Wat} 
\begin{equation}\label{eq:length_parte_real_2}
(a+d)_{\mathbb{R}}=\left(\lambda +\lambda^{-1}\right)\prod_{i=1}^{\lfloor\frac{n-3}{2}\rfloor} \cos(\theta_i),
\end{equation}
where $\lambda$ is the multiplier of $r$, and $2\theta_1$, $2\theta_2, \ldots, 2\theta_{\lfloor\frac{n-3}{2}\rfloor}$ denote its rotational angles. In this case, $\lambda=e^{\frac{\ell(r)}{2}}$ and we get by \eqref{eq:length_parte_real_1} and \eqref{eq:length_parte_real_2} that

$$r_{\mathbb{R}}=\cosh\left(\frac{l(r)}{2}\right)\prod_{i=1}^{\lfloor\frac{n-3}{2}\rfloor} \cos(\theta_i).$$
Hence, $|r_{\mathbb{R}}|\leq \cosh\left(\frac{l(r)}{2}\right)$, and the proposition follows.
\end{proof}

%% file: Sec4_5_Congruence_subgroups.tex
\subsection{Arithmetic Hyperbolic Manifolds} 
We recall that a discrete subgroup $\Gamma \subset \Spin(1,n)$ is {\it arithmetic} if there exist a number field $k$, a $k$-algebraic group $\mathrm{\textbf{H}}$, and an epimorphism $\varphi:\mathrm{\textbf{H}}(k\otimes_{\mathbb{Q}}\mathbb{R})\rightarrow \Spin(1,n)$ with compact kernel such that $\varphi(\mathrm{\textbf{H}}(\mathcal{O}_k))$ is commensurable to $\Gamma$. Here, $\mathcal{O}_k$ denotes the ring of integers of $k$ and $\mathrm{\textbf{H}}(\mathcal{O}_k)=\mathrm{\textbf{H}}\cap \mathrm{GL}_{n}(\mathcal{O}_k)$  with respect to some fixed embedding of $\mathrm{\textbf{H}}$ into $\mathrm{GL}_{n}$. The field $k$ is the \textit{field of definition} of $\Gamma$. Any algebraic $k$-group $\mathrm{\textbf{H}}$ satisfying these properties is called \textit{admissible}. A hyperbolic orbifold $M=\Gamma\backslash\mathbb{H}^n$ such that $\Gamma$ is an arithmetic
subgroup of $\mathrm{Isom}^{+}(\mathbb{H}^n)$ 
is called an {\it arithmetic hyperbolic orbifold}. 
It follows from Borel and Harish–Chandra's Theorem that any arithmetic hyperbolic orbifold has finite volume \cite{BHCh62}.

\subsection{Arithmetic groups of the first type}\label{arithmetic_groups}
The admissibility condition implies that $k$ is a totally real number field, $\mathrm{\textbf{H}}$ is a simply-connected algebraic $k$-group, and by fixing an embedding $k\subset\mathbb{R}$ we can assume that $\mathrm{\textbf{H}}(\mathbb{R})=\Spin(n,1)$ (see \cite[Sec. 2.2]{BE2012} and \cite[Sec. 13.1]{Eme09}). Suppose $k$ is a totally real number field, and $f$ is a quadratic form over $k$. The admissibility of the algebraic $k$-group $\mathrm{\textbf{H}}=\Spin_{f}$ is equivalent to the fact that $f$ has signature $(n,1)$ over $\mathbb{R}$, and $f^{\sigma}$ is definite for any non-trivial embedding $\sigma:k\rightarrow\mathbb{R}$, where \(f^{\sigma}\) is the quadratic form defined on the same \(k\)-space of \(f\), given by \(f^{\sigma}(x)=\sigma(f(x))\). A quadratic form satisfying these conditions will be called an \textit{admissible quadratic form}. The arithmetic groups commensurable to $\Spin_{f}(\ok)$ with $f$ admissible, are called \textit{arithmetic groups of the first type}. For $n$ even, any arithmetic subgroup of $\Spin(n,1)$ is of the first type. There is a second class in odd dimensions, arising from skew-hermitian forms over division quaternion algebras. For $n=7$ there is a third class, arising from certain Cayley algebras. In this work we will deal with arithmetic hyperbolic orbifolds of the first type in dimensions $n\geq 2$.

\subsection{Congruence subgroups}
Let $\Gamma$ be an arithmetic subgroup of $\Spin(n,1)$ commensurable with $\varphi(\mathrm{\textbf{H}}(\mathcal{O}_k))$, and $M=\Gamma\backslash\mathbb{H}^n$ the corres\-ponding hyperbolic arithmetic orbifold.
If $I\subset\mathcal{O}_{k}$ is a non-zero
ideal of $\mathcal{O}_{k}$, the {\it principal congruence subgroup} of $\Gamma$ associated to $I$ is the subgroup $\Gamma(I):=\Gamma \cap \varphi(\mathrm{\textbf{H}}(I))$, where 
$$\mathrm{\textbf{H}}(I):=\mathrm{ker}\left(\mathrm{\textbf{H}}(\mathcal{O}_{k})\xrightarrow{\pi_{I}} \mathrm{\textbf{H}}(\mathcal{O}_{k}/I)\right),$$
and $\pi_{I}$ denotes the reduction modulo $I$ map. The associated {\it principal congruence covering} is $M_I=\Gamma(I)\backslash\mathbb{H}^n \rightarrow M$. Since $\Gamma(I)$  is a normal finite-index subgroup of $\Gamma$, the covering $M_I\rightarrow M$ is a regular finite-sheeted covering map. More generally, a discrete subgroup $\Lambda < \Gamma$ is called a {\it congruence subgroup} if $\Gamma(I)\subset\Lambda$ for some ideal $I\subset \mathcal{O}_{k}$. 

Let $f$ be an admissible quadratic form over a totally real number field $k$ of degree $d$. We can describe the group $\Gamma=\Spin_f(\ok)$ and its principal congruence subgroups $\Gamma(I)$ in the following way. Denote by $e_1, e_2,\ldots,e_{n+1}$ an orthogonal basis with respect to $f$. Then under the linear representation given by left multiplication in $\mathscr{C}^{+}(f,\mathbb{R})$ we get

$$\Gamma=\left\lbrace s=\sum_{\mathclap{|M| \mbox{ even }}}s_{M}e_{M} ~ | ~  s_{M} \in\ok \mbox{ and } ss^*=1 \right\rbrace$$
and
$$\Gamma(I)=\left\lbrace s=\sum_{\mathclap{|M| \mbox{ even}}} s_{M}e_{M}\in \Gamma ~ | ~ s_{M}\in I \mbox{ for } M\neq\emptyset \mbox{ and } s_{\mathbb{R}}-1\in I \right\rbrace.$$
(see \cite[Sec. 2.4]{Mur19}). To simplify the discussion, we will denote by $\mathcal{Q}$ the $\ok$-order in $\mathscr{C}^{+}(f,\mathbb{R})$ given by
$$\mathcal{Q}=\left\{s=\sum_{\mathclap{|M| \mbox{ even}}}s_{M}e_{M} ~ | ~ s_{M}\in\ok\right\}.$$

\subsection{Length inequality for $\Gamma(\alpha)$}
For a principal ideal $I=(\alpha)$, $\alpha\in\ok$, we denote $\Gamma(I)$ simply by $\Gamma(\alpha)$. In the same way, we denote by \(\N(\alpha)\) the \emph{norm} of the ideal \((\alpha)\). We will present a series of results relating the real part of elements in $\Gamma(\alpha)$ with $\alpha$ that will be necessary in the course of the investigation.

\begin{lemma}\label{realpart}
Let $\alpha \in \ok$ be a nonzero element. For any $s \in \Gamma(\alpha)$ we have the equality
\[ s_\mathbb{R}=1+\frac{1}{2}\alpha^2 \zeta \]
for some $\zeta \in \ok.$
\end{lemma}
\begin{proof}
 By definition, we can write $s=1+\alpha t$ for some $t \in \mathcal{Q}$. Since $s^*=1+\alpha t^*$ we have
 \[1=s s^* =1+\alpha(t+t^*)+\alpha^2 tt^*.\]
Now, $t\in\mathcal{Q}$ implies that $(tt^*)_{\mathbb{R}}$ lies in $\mathcal{O}_k$. By taking the equality of real parts, and observing that $2t_\mathbb{R}=(t+t^*)_\mathbb{R}$, the lemma follows with $\zeta=-(tt^*)_\mathbb{R}$.
\end{proof}
 
The following complement of Lemma \ref{realpart} will also be useful.

\begin{lemma}\label{complement_real_part}
 Let $\alpha \in \ok$ be a  nonzero element and $s \in \Gamma$  such that $s-s_\mathbb{R} \in \alpha \mathcal{Q}$. For any $r \in \Gamma(\alpha)$ we have the equality
 \[ (sr)_\mathbb{R}=s_\mathbb{R}+\frac{1}{2}\alpha^2 \zeta\]
 for some $\zeta \in \ok.$
\end{lemma}
\begin{proof}
 Again, we have $r=1+\alpha t$ for some $t \in \mathcal{Q}$ with $2t_\mathbb{R}=\alpha \xi$ for some $\xi \in \ok.$ Since $s=s_\mathbb{R}+\alpha u$ for some $u \in \mathcal{Q}$ we can write 
 \[sr=s(1+\alpha t)=s+\alpha(s_\mathbb{R}+\alpha u)t=s+\alpha s_\mathbb{R} t+\alpha^2ut.\]
Since $ut\in\mathcal{Q}$, then $(ut)_{\mathbb{R}}$ lies in $\mathcal{O}_k$. When we take the real part in the last equality we finish the proof with $\zeta=s_\mathbb{R} \xi+2(ut)_\mathbb{R}.$
\end{proof}

\begin{prop}\label{lower_bound_for_real_part}
  Let $\alpha \in \ok$ be a nonzero element. For any loxodromic element $r \in \Gamma(\alpha)$ we have
 \[|r_\mathbb{R}| > \frac{1}{2^{2d-1}}\mathrm{N}(\alpha)^2-1,\]
where $d=[k:\Q]$.
\end{prop}
\begin{proof}
See \cite[Lemma 4.1]{Mur19}.
\end{proof}

\subsection{The congruence subgroup $\Gamma_\tau(\alpha)$}
An important example of congruence subgroup which will play a special role in this work is the following. Fix an element $\alpha\in \ok$. Suppose that  $\tau \in (\ok/\alpha\ok)^\times$ is an element of order 2, and define
$$\Gamma_\tau(\alpha)=\{\gamma \in \Gamma \mid \gamma \in \Gamma(\alpha)\hspace{2mm} \mbox{or} \hspace{2mm}\gamma\equiv \tau (\mathrm {mod}~ \alpha \mathcal{Q}) \}.$$

The group $\Gamma_\tau(\alpha)$ is a normal subgroup of $\Gamma$. Indeed, let \(t \in \ok\) such that \(\tau=t+\alpha \ok\). We note that $\{1 + \alpha \mathcal{Q},t+\alpha \mathcal{Q}\}$ is a \emph{central} subgroup of \(\left(\mathcal{Q}/\alpha\mathcal{Q}\right)^{\times}\), hence $\Gamma_\tau(\alpha)$ is normal since it is the preimage of this normal subgroup under the natural projection map $\Gamma\to\left(\mathcal{Q}/\alpha\mathcal{Q}\right)^{\times}$ .
 
\section{Hyperbolic manifolds with a systole lying in a surface}
The goal of this section is to show that, under certain conditions, the manifold $\Gamma_{\tau}(\alpha)\backslash\mathbb{H}^{n}$ has a systole contained in a totally geodesic surface. 

\subsection{A totally geodesic surface embedded in $\Gamma_{\tau}(\alpha)\backslash\mathbb{H}^{n}$}\label{sec:surface_embedding}

Let $k$ be a totally real number field, and $(E,f)$  be an admissible $(n+1)$-dimensional quadratic space over $k$ (see subsection 4.2). Since \(f\) has signature \((n,1)\) and \(k \subset \R\), by the Gram-Schmidt process and the Law of Inertia, there exists a basis $\{e_0,e_1,\ldots,e_n\}$ of \(E\) such that in this basis $f$ has the diagonal form $f=-a_0x_0^2+a_1x_1^2+\cdots + a_nx_n^2$, with $a_i \in k$ and positive for all \(i\). In fact, we can suppose that the coefficients are in \(\ok\) if we replace \(e_i\) by \(d_ie_i\) where \(d_i \in \ok\) is a denominator of \(a_i\).
By admissibility, \(f^\sigma\) is positive definite for any non-trivial Galois embedding $\sigma:k\rightarrow \mathbb{R}$, hence $\sigma(a_0)=-\sigma(-a_0)=-f^{\sigma}(e_0)<0$ and $\sigma(a_i)=f^{\sigma}(e_i)>0$ for all $i=1,\ldots,n$. 

 Let $E'$ be the subspace generated by $\{e_0, e_1,e_2\}$, and $f':E'\rightarrow k$ the restriction of $f$ to $E'$. The inclusion $E'\rightarrow E$ defines a natural inclusion $\Gamma'=\Spin_{f'}(\ok)\hookrightarrow\Spin_f(\ok)=\Gamma$. For any $\alpha\in\ok$ and $\tau\in(\ok\backslash\alpha\ok)^{\times}$ of order two, by definition we get an inclusion $$\Gamma'_{\tau}(\alpha)\hookrightarrow\Gamma_{\tau}(\alpha).$$ Consider an isometric embedding of $\mathbb{H}^2$ into $\mathbb{H}^n$ equivariant with res\-pect to the actions of $\Gamma'$ and $\Gamma$ and the inclusions above. For any $\alpha$ and $\tau$ as before, we obtain a totally geodesic embedding 

\begin{equation}\label{eq:surface_embedding}
S_{\alpha,\tau}\hookrightarrow M_{\alpha,\tau},
\end{equation} 
where $S_{\alpha,\tau}=\Gamma'_{\tau}(\alpha)\backslash\mathbb{H}^2$ and $M_{\alpha,\tau}=\Gamma_{\tau}(\alpha)\backslash\mathbb{H}^n$. This implies in parti\-cular that $\sys(M_{\alpha,\tau})\leq\sys(S_{\alpha,\tau})$.

\begin{prop}\label{complement_lower_bound_for_real_part}
	Let $\alpha \in \ok$ be a  nonzero element and $s \in \Gamma'$  such that $s-s_\mathbb{R} \in \alpha \mathcal{Q} $. Then $\tau=\overline{s_\mathbb{R}}$ has order two in $(\ok/\alpha\ok)^\times$. Furthermore, for any loxodromic element $\gamma \in \Gamma_\tau(\alpha) \setminus \Gamma(\alpha)$ we have  
	\[|\gamma_\mathbb{R}| > \frac{1}{2^{2d-1}}\mathrm{N}(\alpha)^2-|s_\mathbb{R}|.\]
\end{prop}
\begin{proof}
Indeed, since \(s\) is contained in a quaternion algebra we have $s=s_\mathbb{R}+\alpha u$ for some $u \in \mathcal{Q}$ with $u^*+u=0$. Hence, the equation $1=ss^*=s_\mathbb{R}^2+\alpha^2uu^*$ implies that $s_\mathbb{R}^2=1 (\mathrm{mod}~ \alpha)$.
	Since the index $[\Gamma_\tau(\alpha): \Gamma(\alpha)]=2$, we need to estimate the real part of any product $\gamma=sr$ with $r \in \Gamma(\alpha)$. In this case, by Lemma \ref{complement_real_part}, we get 
	\begin{equation}\label{eq:partes reais}
	\gamma_\mathbb{R}=s_\mathbb{R}+\frac{1}{2}\alpha^2\zeta.
	\end{equation}
	
	Now, for any non-trivial archimedean place $\sigma$ of $k$ we know that $|\sigma(\gamma_\mathbb{R})|\leq 1$, and $|\sigma(s_\mathbb{R})|\leq 1$ (\cite[Equation 8]{Mur19}). Therefore, by applying $\sigma$ to Equation \eqref{eq:partes reais} we get
	\begin{equation*}
	|\sigma(\alpha^2\zeta)|=2|\sigma(\gamma_\mathbb{R})-\sigma(s_{\mathbb{R}})|\leq 4 .
	\end{equation*}
	Once more, by \eqref{eq:partes reais}, and the fact that $\zeta\in\ok$, we obtain that
	\begin{align*}
	|\gamma_{\mathbb{R}}|&\geq \frac{1}{2}|\alpha|^2|\zeta| -|s_{\mathbb{R}}|\\
	&=\frac{1}{2\prod_{\sigma\neq id}|\sigma(\alpha)^2\sigma(\zeta)|}\N(\alpha)^2|\N(\zeta)|-|s_{\mathbb{R}}|\\
	&\geq\frac{1}{2^{2d-1}}\N(\alpha)^2 - |s_{\mathbb{R}}|.
	\end{align*}
\end{proof}

\subsection{Embeddings of quadratic fields in $\mathscr{C}^{+}(f',k)$}
A direct computation shows that the Clifford algebra $\mathscr{C} ^{+}(f',k)$ is a quaternion algebra of $\Gamma'$ (see \cite[Section 12.2]{MR}). In fact, it coincides with the invariant quaternion algebra. It is well known that closed geodesics in $S'=\Gamma'\backslash\mathbb{H}^2$ are related with quadratic extensions of $k$ that embed in $\mathscr{C} ^{+}(f',k)$. In this subsection we will recall the important properties of this connection that will be useful in the sequel.

For any $s \in \mathscr{C} ^{+}(f',k)$ the \textit{reduced trace} $a=s+s^*$, and the \textit{reduced norm} $b=ss^*$ are elements of $k$. Hence $s$ is a root of the quadratic polynomial $g(x)= x^2-ax+b \in k[x]$. If $g$ is irreducible over $k$, and $L \supset k$ is the quadratic extension where $g$ splits, then for any fixed root $\alpha$ of $g$ in $L$, there exists a unique monomorphism $\phi:L \to \mathscr{C} ^{+}(f',k)$ such that $\phi(\alpha)=s$, $\phi\mid_{k}$ is the identity, and $\phi(\sigma(x))=\phi(x)^*$ for the non-trivial Galois automorphism $\sigma:L\rightarrow L$ of $L$ over $k$. 
In particular, via the identification of $L$ with $\phi(L)$, the map $\sigma$ coincides with the restriction of $^*$ to $L$.

The following proposition is a well-known fact about quaternion algebras which we recall here for reader convenience.

\begin{prop}\label{quadratic_extension}
Let $s\in\Gamma'$ be a loxodromic element. There exist a quadratic extension $L=k(\sqrt{D})$ for some $D \in k$ positive, and a $k$-homomorphism $\psi:L \to \mathscr{C} ^{+}(f')$ such that $s=s_\mathbb{R}+\psi(\sqrt{D}).$
\end{prop}
\begin{proof}
Consider the irreducible polynomial 
\begin{equation}{\label{eq:minimal-polynomial-dim2}}
g(x)=x^2-(s+s^*)x+1
\end{equation}
over $k$. Since $s$ is loxodromic, $g$ has two distinct real roots $\lambda$ and $\lambda^{-1}$, and let $L$ be the quadratic extension $k(\lambda)$. Without loss of generality, suppose that $|\lambda|>1$. 
We note that $\lambda$ and $s$ are roots of $g$, thus there is a unique $k$-homomorphism $\phi:L \to \mathscr{C} ^{+}(f')$ with $\phi(\lambda)=s$. If we define $\theta=\lambda-s_\mathbb{R} \in L$, the equality $\lambda+\sigma(\lambda)=2s_\mathbb{R}$ implies $\theta+\theta^*=0$, and then $\theta^2=s_{\mathbb{R}}^{2}-1=:D$. Moreover, $s$ loxodromic implies that $|2s_{\mathbb{R}}|>2$, thus $D>0$.
To finish the proof, if $\theta>0$ we take $\psi=\phi$, otherwise we consider $\psi=\phi^*$.
\end{proof}

By Proposition \ref{quadratic_extension}, we can write $s=\psi(\lambda_0)$ where $\lambda_0=s_\mathbb{R}+\sqrt{D}$. Thus we have an isomorphism between the cyclic group generated by $\lambda_0$ in $L=k(\sqrt{D})$, and the cyclic group generated by $s$ in $\Gamma'$. For each $n \in \mathbb{N}$ we can write $\lambda_0^{n+1}=t_n+u_n\sqrt{D}$ with $t_n,u_n \in \ok$. In the next result we obtain asymptotic relations between $u_n, t_n$ and $s_{\mathbb{R}}$.

\begin{lemma}\label{asymptotic_of_t_and_u}
 If $\lambda=x_0+\sqrt{D}$ is a unit in $\ok[\sqrt{D}]$ and $\lambda^{n+1}=t_n+u_n\sqrt{D},$ then for each $n \geq 1$, we have
 \begin{equation}\label{asymptotic_tn}
  t_n=2^nx_0^{n+1}+O(x_0^n)
 \end{equation}
and
\begin{equation}\label{asymptotic_un}
 u_n=2^nx_0^n+O(x_0^{n-1}),
\end{equation}
where the O notation is considered with respect to \(x_0.\)
\end{lemma}
\begin{proof}
 Since $D=x_0^2-1$, we have the following relations
 \[t_{n}=(x_0^2-1)u_{n-1}+x_0t_{n-1} \mbox{ and } u_{n}=x_0u_{n-1}+t_{n-1}.\]
 Hence,
 \[t_n=x_0u_{n}-u_{n-1} \mbox{ and } u_n=2x_0u_{n-1}-u_{n-2} \mbox{ for all } n \geq 2.\]
 We prove \eqref{asymptotic_un} by induction. For $n=0, 1$ the relation is trivial. Assuming valid for any $1 \leq k <n$, it follows that
 \[u_n=2x_0(2^{n-1}x_0^{n-1}+O(x_0^{n-2}))-2^{n-2}x_0^{n-2}+O(x_0^{n-3})=2^nx_0^n+O(x_0^{n-1}).\]
Now, we obtain \eqref{asymptotic_tn} from the relation $t_n=x_0u_{n}-u_{n-1}$.
\end{proof}

Recall that a real algebraic integer $\lambda >1$ is a \textit{Salem number} if $\lambda^{-1}$ is a Galois conjugate of $\lambda$, and the other conjugates of $\lambda$ lie on the unit circle. It is well known that the roots of the characteristic polynomial associated to loxodromic elements in $\Gamma'$ are \textit{Salem numbers} (see \cite{GH01}).  In this work, we will only deal with Salem numbers of degree four. The next subsection contains the main result that will be necessary for our purpose.  

\subsection{Salem numbers of degree four}
For the interest of this work, it is important to develop some results about Salem numbers of low degree.  
Let $\mu$ be a Salem number of degree four. The field $K=\mathbb{Q}(\mu+\mu^{-1})$ is a totally real number subfield of $\mathbb{Q}(\mu)$,  with nontrivial $\mathbb{Q}$-isomorphism $\tau:K \to K$. Since $[\mathbb{Q}(\mu):K]=2$, there exists a unique nontrivial $K$-isomorphism $\sigma:\mathbb{Q}(\mu) \to \mathbb{Q}(\mu)$ such that $\sigma(\mu)=\mu^{-1}$. Hence, the four embeddings of $\mathbb{Q}(\mu)$ into $\mathbb{C}$ are the inclusion, $\sigma$, $\tau$ and $\overline{\tau}$, where $\tau$ is the extension of the nontrivial $\mathbb{Q}$-morphism of $K$ into $\mathbb{Q}(\mu)$. In particular, we can assume that $\tau(\mu)=e^{i\nu}$ for some $\nu \in (0,\pi)$.

Now, suppose that there exists $D \in \mathcal{O}_K$ such that $\mathbb{Q}(\mu)=K(\sqrt{D})$ and $\mu=t+u\sqrt{D}$ for some $t,u \in \mathcal{O}_K$. Since $\sigma$ is the non-trivial $K$-automorphism of $\mathbb{Q}(\mu)$ then $\sigma(\sqrt{D})=-\sqrt{D}$, and
\begin{equation}
1=\mu\cdot \mu^{-1}=\mu \cdot \sigma(\mu)=(t+u\sqrt{D})\cdot(t-\sqrt{D})=t^2 - u^2D
\end{equation}
\begin{equation}\label{angle_conjugated}
2\tau(t)=\tau(2t)=\tau(\mu+\sigma(\mu))=\tau(\mu)+\tau(\mu^{-1})=2\cos(\nu).
\end{equation}

Thus $t^2-u^2D=1$ and $|\tau(t)| < 1$. For geometric reasons, it is important to get Salem numbers in this form such that $\tau(t)$ is not very small. The next proposition shows that we can assume that this property is true up to a small power of $\mu$.

\begin{prop}\label{Salem_angles}
Let $\mu>1$ be a Salem number of degree four. With the previous notations, if $\mu=t+u\sqrt{D}$, $t,u\in\mathcal{O}_{K}$ there exists $m \in \{0,1,2\}$ such that $\mu^{m+1}=t_m+u_m\sqrt{D}$ with $\tau(t_m)^2>\frac{1}{2}.$

\end{prop}
\begin{proof}
 For each $m \in \{0, 1, 2\}$, if $\mu^{m+1}=t_m+u_m\sqrt{D}$ with $t_m,u_m \in \mathcal{O}_K$, then \[\mu^{2(m+1)}=(t_m^2+Du_m^2)+2t_mu_m\sqrt{D}=(2t_m^2-1)+2t_mu_m\sqrt{D}.\]

Hence, $t_{2m+1}=2t_m^2-1$ and $\tau(t_{2m+1})=2\tau(t_m)^2-1$. Hence $\tau(t_m)^2>\frac{1}{2}$ if, and only if, $\tau(t_{2m+1})>0.$ On the other hand,
\[\tau(2t_{2m+1})=\tau(\mu^{2(m+1)})+\mu^{-2(m+1)})=2\cos(2(m+1)\nu),\]
where $m+1 \in\{1,2,3\}$.
Then, it remains to show that $\cos(2k\nu)>0$ for some $k \in \{1,2,3\}$. Indeed, consider the sets $S_1=(0,\frac{\pi}{4}) \cup (\frac{3\pi}{4},\pi), S_2=(\frac{3\pi}{8},\frac{5\pi}{8}), S_3=(\frac{\pi}{4},\frac{3\pi}{8}) \cup (\frac{5\pi}{8},\frac{3\pi}{4})$. 
For each $\nu \in S_j$, we have $\cos(2j\nu)>0$. The lemma is now proven since $[0,\pi]-(S_1 \cup S_2 \cup S_3)$ only contains rational multiples of $\pi$ and Salem numbers do not have conjugates of finite order.
\end{proof}

\subsection{Congruence subgroups with explicit elements of minimal displacement}
The purpose of this section is to construct hyperbolic manifolds with systole lying in a totally geodesic surface. More specifically, we are looking for conditions on $\tau$ and $\alpha$ such that the manifold $M_{\tau,\alpha}$ has a systole in $S_{\tau,\alpha}$ (see Section \ref{sec:surface_embedding}). Since $\sys(M_{\alpha,\tau})\leq\sys(S_{\alpha,\tau})$, it is necessary to bound $\sys(M_{\alpha,\tau})$ from below. Proposition \ref{lower_bound_for_real_part}  and Proposition \ref{complement_lower_bound_for_real_part}  show that this require a lower bound for the norm of $\alpha$ in the base field $k$, which at the same time implies that the Galois conjugates of $\alpha$ cannot be very small. We are able to find such $\alpha$ when $k$ is a \textit{real quadratic number field}.

In the sequel, we consider the definitions of $\Gamma$, $\Gamma'$ as in Section \ref{sec:surface_embedding}, and $k$ a real quadratic field, with $\sigma:k\rightarrow\mathbb{R}$ the nontrivial embedding of $k$ into $\mathbb{R}$.   

\begin{lemma}\label{level_lemma}
Let $s \in \Gamma'$ be a primitive loxodromic element with $s_\mathbb{R}>0$. There exists $l \in \{2,5,8\}$ which depends only on $s_\mathbb{R}$ such that $s^{l+1}=(s^{l+1})_\mathbb{R}+\alpha_l u_l$ with $u_l \in \mathcal{Q}$, $\alpha_l\in \ok$ and
\[1\leq |\sigma(\alpha_l)| \leq 5.\] Moreover, if $t=s_\mathbb{R}$, then the following asymptotic relation
 \[\alpha_l = C_lt^{\frac{2}{3}(l+1)}+O\left(t^{\frac{2}{3}(l+1)-2}\right),\]
holds for some constant $C_l>0$ which depends only on $l$.
\end{lemma}
\begin{proof}
By Proposition \ref{quadratic_extension} and the discussion in Section \ref{sec:surface_embedding}, the element $s$ corresponds to a Salem number $\lambda_0=t_0+\sqrt{D}$, and $L=k(\lambda_{0})$ is a quadratic extension of $k$.
Since $k$ is a real quadratic number field, $\lambda_0$ is a Salem number of degree four. By Proposition \ref{Salem_angles}, there exists $m \in \{0,1,2\}$ such that $\lambda_0^{m+1}=t_m+u_m\sqrt{D}$ with $|\sigma(t_m)|^2>\frac{1}{2}$. For convenience, we can rewrite \[\lambda=\lambda_0^{m+1}=t+\sqrt{E}\]
 where $t=t_m$ and $E=u_m^2D$. 
It is straightforward to check that
 \[\lambda^3=(4t^3-3t)+(4t^2-1)\sqrt{E}.\]
 
If $l$ is given by $l=3(m+1)-1,$ then
 \[\lambda_0^{l+1}=t_l+u_l\sqrt{D}=t_l+\alpha_l\sqrt{E}\]
 where $\alpha_l=4t_m^2-1$ and $E=u_m^2D$. Since $\frac{1}{2}<|\sigma(t_m)|^2<1$ we conclude that $1 < |\sigma(\alpha_l)| < 5$. The asymptotic behaviour of $\alpha_l$ follows directly from \eqref{asymptotic_tn} and the equality $m+1=\frac{1}{3}(l+1).$
\end{proof}

We will now prove that for any primitive element $s \in \Gamma'$ producing a closed geodesic with length sufficiently large, some power $s^l$ with $l$ \textit{uniformly bounded} realizes the systole of some congruence hyperbolic $n$-manifold. 

\begin{prop}\label{inedisp}
There exists a universal constant $L>0$ such that for any loxodromic element $s\in \Gamma'$ with $s_{\mathbb{R}}>L,$ we can find $l \in \{2,5,8\}$ depending only on $s_\R$ with $s^{l+1}-(s^{l+1})_\mathbb{R} \in \alpha_l\mathcal{Q}$, for some \(\alpha_l \in \ok\) and
	\[\ell(s^{l+1})\leq \ell(r) \mbox{ for all loxodromic element } r \in \Gamma_{\tau_l}(\alpha_l),\]
	where $\tau_l$ is the class of $(s^{l+1})_\mathbb{R}$ modulo $\alpha_l$.
\end{prop}

\begin{proof}
Fix $s \in \Gamma'$ loxodromic with real part $s_\R$. By Lemma \ref{level_lemma} there exists $l \in \{2,5,8\}$ depending only on $s_\R$ such that $s^{l+1}-(s^{l+1})_\mathbb{R} \in \alpha_l\mathcal{Q}$ for some \(\alpha_l \in \ok\) with \(1 \leq |\sigma(\alpha_l)| \leq 5\). 
 Let $t_l=(s^{l+1})_\mathbb{R}$, it follows from Proposition \ref{complement_lower_bound_for_real_part} that $t_l^2 \equiv 1 (\mathrm{mod}~ \alpha_l)$. Hence, if we denote by $\tau_l$ the class of $t_l$, we have $$\Gamma_{\tau_l}(\alpha_l)=\Gamma(\alpha_l) \cup s^{l+1}\Gamma(\alpha_l).$$ 

If $r \in \Gamma(\alpha_l)$, since $|\alpha_l'| \geq 1$, Proposition \ref{lower_bound_for_real_part} gives us that
	\[ |r_\mathbb{R}| \geq \frac{1}{8}|\sigma(\alpha_l)|^2\alpha_l^2-1 \geq \frac{1}{8}\alpha_l^2-1.\]

Since $\ell(s^{l+1})=2\cosh^{-1}(|t_l|)$, by Proposition \ref{prop:realpart-length}, in order to show that $\ell(r) \geq \ell(s^{l+1})$ it is sufficient to guarantee that $|r_\mathbb{R}| \geq t_l$. 
	By Lemma \ref{level_lemma} and Equation \eqref{asymptotic_tn} we have
	\[\frac{1}{8}\alpha_l^2-1 = C(s_\R)^{\frac{4}{3}(l+1)} + O\left((s_\R)^{\frac{4}{3}(l+1)-4}\right) \mbox{ and } t_l \sim 2^l(s_\R)^{l+1}. \]
	Hence, $\frac{1}{8}\alpha_l^2-1 \geq t_l$ whenever $s_\R$ is sufficiently large.

Analogously, if $r \in s^{l+1}\Gamma(\alpha_l)$ and $r \neq s^{l+1}$, we have by Lemma \ref{complement_real_part} that
\[ |r_\mathbb{R}| \geq \frac{1}{8}|\sigma(\alpha_l)|^2\alpha_l^2-t_l \geq \frac{1}{8}\alpha_l^2-t_l.\]
	And $|r_\mathbb{R}|>t_l$ whenever $\alpha_l^2>16t_l,$ which holds whenever $s_\R$ is large enough.
\end{proof}

The previous result has Theorem \ref{th:embeded_surface} mentioned in Introduction as geometric counterpart. We recall it as a corollary, adapted to the terminology used so far.

\begin{corollary}\label{equalsystoles}
With the notation as in Proposition \ref{inedisp}, the  hyperbolic manifold $M_{\tau_{l},\alpha_{l}}=\Gamma_{\tau_{l}}(\alpha_{l})\backslash\mathbb{H}^{n}$ contains the totally geodesic surface $S_{\tau_{l},\alpha_{l}}=\Gamma'_{\tau_{l}}(\alpha_{l})\backslash\mathbb{H}^{2}$ with $$\sys(S_{\tau_{l},\alpha_{l}})=\sys(M_{\tau_{l},\alpha_{l}}).$$ 
\end{corollary}
\begin{proof}
Indeed, by Proposition \ref{inedisp}, we exhibit a loxodromic element in \(\Gamma'_{\tau_l}(\alpha_l)\), namely \(s^{l+1}\), of minimal displacement in \(\Gamma_{\tau_l}(\alpha_l)\). In particular
this element is also of minimal displacement in \(\Gamma'_{\tau_l}(\alpha_l)\). Hence, its length is at same time the systolic length of \(S_{\tau_{l},\alpha_{l}}\) and \(M_{\tau_{l},\alpha_{l}}.\)
\end{proof}

%% file: Sec6_Proofs.tex
Any closed geodesic on a hyperbolic manifold \(M\) is parametrized by a constant speed curve from the circle to \(M\) and we can identify the geodesic with the equivalence class of such parametrization up to re\-parametrization. 
Let $\gamma:\mathbf{S}^1=\mathbb{R}/\Z \rightarrow M$ be a closed geodesic, we say that $\gamma$ is \emph{primitive} if \(\gamma\) is injective, i.e. if \(\gamma\) is an embedding. Any closed geodesic $\delta$ is a $k$-fold iterate of some primitive geodesic $\gamma$, i.e. there exists $k\in \mathbb{N}$ such that $\delta(t)=\gamma(kt)$ (up to reparametrizations of \(\delta\) and \(\gamma\)). We note that \(k\) is uniquely determined by the relation \(\ell(\delta)=k\ell(\gamma)\) and because of this, we call it as the {\it order} of $\delta$.

Let $\pi: M \to N$ be a covering map  between two hypebolic $n$-orbifolds $M$ and $N$, we recall that a closed geodesic $c: \mathbf{S}^1 \to N$ {\it lifts} to $M$ if there is a closed geodesic $\tilde{c}: \mathbf{S}^1 \to M$ such that \(c=\pi \circ \tilde{c}\). In this case, we say that \(\tilde{c}\) is a lift of $c$. \\
\indent We note that the deck group $\mathrm{Deck}(\pi)=\{ g \in \Isom(M) \mid \pi \circ g = \pi \}$ is always finite whenever $M$ and $N$ have finite volume. In the sequel, we consider the natural action of $\mathrm{Deck(\pi)}$ on the set of closed geodesics of $M$.

\begin{lemma}\label{C051121}
	Let $\pi: M \to N$ be a covering map between the hyperbolic $n$-orbifolds $M$ and $N$ of finite volume, and let $G=\mathrm{Deck}(\pi)$.
	\begin{enumerate}
		\item If $\tilde{\gamma}_1,\tilde{\gamma_2}$ are closed geodesics on $M$ which are liftings of two distinct closed geodesics $\gamma_1,\gamma_2$ on \(N\) respectively, then the orbits $G \cdot \tilde{\gamma}_1$ and $G \cdot \tilde{\gamma}_2$ are disjoint.
		\item If $\gamma$ is a closed geodesic on \(N\) of order $k$ which lifts, then for any lifting $\tilde{\gamma}$ its isotropy group $G_{\tilde{\gamma}}$ has at most $k$ elements. 
	\end{enumerate} 
\end{lemma}
\begin{proof}
	Indeed, if $\tilde{\gamma}_1 = g \circ \tilde{\gamma}_2$ (up to reparametrization of \(\tilde{\gamma_1}\), and \(\tilde{\gamma_2}\)) then $\gamma_1=\pi \circ \tilde{\gamma}_1=\pi \circ g\circ \tilde{\gamma}_2=\pi \circ \tilde{\gamma}_2=\gamma_2,$ which proves $(1)$. For $(2)$ we can suppose that $M=\Lambda' \backslash \HS^n$ and $N=\Lambda \backslash \HS^n$ where $\Lambda' < \Lambda$.  With this identification, the group $G$ can be seen as $\mathrm{N}_\Lambda(\Lambda')/\Lambda'$, where \(\mathrm{N}_\Lambda(\Lambda')\) denotes the \emph{normalizer} of \(\Lambda'\) in \(\Lambda\). Moreover, a closed geodesic on $M$ can be associated  with a conjugacy class $[\gamma']$ of a loxodromic element $\gamma' \in \Lambda'$. The action of $G$ on the set of closed geodesics is given by $\lambda\Lambda' \cdot [\gamma']= [\lambda\inv \gamma' \lambda]$. If $[\gamma']$ denotes a closed geodesic of order $n$ on $N$, we can use the same notation for its lifting on $M$ since $\gamma \in \Lambda'$. Hence $\lambda \Lambda' \cdot [\gamma'] =   [\gamma']$ means that  $\lambda \inv \gamma' \lambda=\lambda_1\inv \gamma' \lambda_1$ for some $\lambda_1 \in \Lambda'$, then $\lambda_1\lambda\inv$ commutes with $\gamma'$. By hypothesis, $\gamma'=\eta_0^n$,  and for results in hyperbolic geometry we have that the centralizer of $\gamma'$ is the cyclic group generated by $\eta_0$. Therefore, $\lambda\Lambda' \in \{\eta_0^i\Lambda' \mid 0 \leq i \leq k-1\}$.   
\end{proof}

\begin{remark}\label{rmk}
		Let $M$ be a closed hyperbolic $n$-manifold and let $\Sigma \subset M$ be a totally geodesic submanifold. If $\alpha,\beta$ are distinct primitive closed geodesics on $\Sigma$ then \(\alpha\) and \(\beta\) are distinct primitive closed geodesics in $M$. Indeed, if $\alpha$ is a $k$-folded iterated of $\alpha_0$ for some  primitive $\alpha_0:\mathbf{S}^1 \to M$, we have $\alpha_0(0) \in \Sigma$ and $\alpha_0'(0) \in T_{\alpha_0(0)}\Sigma$, thus $\alpha_0$ is a closed geodesic on $\Sigma$ and then $\alpha=\alpha_0$. In particular, if $\sys(\Sigma)=\sys(M)$ then $\Kiss(M) \geq \Kiss(\Sigma)$.
\end{remark}

\begin{proof}[Proof of Theorem \ref{firstresult}]
		For each $n \geq 2$ we can consider a fixed closed arithmetic hyperbolic $n$-manifold of the first type $M$. By \cite[Main Theorem]{xue92} there exists a sequence \(M_j\) of congruence coverings of \(M\) such that $\beta_1(M_j) \to \infty$, since \(\vol(M_j) \to \infty\). Moreover, it follows from \cite[Theorem 6.1]{Pli} that $\sys(M_j) \to \infty$. Hence, given $A>0$ we can suppose that we have a closed hyperbolic $n$-manifold $M$ with $\sys(M)>A$ and $\beta_1(M)>0$.
	    
	    If we write $M=\Gamma \backslash \HS^n$ where $\Gamma \simeq \pi_1(M)$, then there exists an epimorphism $\phi:\Gamma \to \Z$. Let $N_t \to M$ be the cyclic cover of degree $t$ obtained by the kernel of the map \[\Gamma \to \Z \to \Z/t \Z.\]                                                                                                                                                                                                        
		
		Let $\gamma \in \Gamma$ be a  nontrivial loxodromic element of minimal displacement satisfying $\phi(\gamma)=0$. If $\alpha$ corresponds to the closed geodesic induced by $\gamma$, then $\alpha$ is lifted for all covering $N_t$, thus $A \leq \sys(N_t) \leq \ell(\alpha)$ for all $t \geq 2$.
		
		For any $t$, there exists  a closed geodesic $\alpha_t$ on $M$ whose lifting $\tilde{\alpha_t}$ in $N_t$ satisfies $\sys(N_t)=\ell(\tilde{\alpha_t})=\ell(\alpha_t)$. Let $G_t \simeq \mathbb{Z}/t\mathbb{Z}$ be the deck group of the covering $N_t \to M$. We claim that the isotropy subgroup of $\tilde{\alpha_t}$ under the action of $G_t$ has cardinality at most $c$ for some constant $c>0$ which does not depend on $t$. Indeed,  $\alpha_t=\beta^{m}$ for some primitive closed geodesic $\beta$. Hence, 
		\[m \cdot \sys(M) \leq m \cdot \ell(\beta)=\ell(\alpha_t)=\sys(N_t) \leq \ell(\alpha).\]
		
		Therefore, by Lemma \ref{C051121} the order of the isotropy group of $\tilde{\alpha_t}$ by the action of the deck group $G_t$ is bounded from above by $c=\frac{\ell(\alpha)}{\sys(M)}>0.$
		
		For any $p$ prime with $p>c$ we conclude that the orbit $G_p \cdot \tilde{\alpha_p}$ has $p$ elements. 
		Putting all the information above together, we have proved that for any $p>c$, $\K(N_p) \geq p$ and $\vol(N_p)=p \cdot \vol(M)$. Thus, the theorem is proved with $B=\ell(\alpha)$ and $C=\vol(M)$. 
\end{proof}

We are now ready to present the proof of Theorem \ref{maintheorem}. 

\begin{proof}[Proof of Theorem \ref{maintheorem}] 
Let $\phi'(x)$ denote the number of conjugacy classes of loxodromic elements in $\Gamma'$ with reduced trace equal to $x$. 
By the Prime Geodesic Theorem  there is a sequence $x_i\rightarrow\infty$ such that (see \cite{Sch95}) 
$$\phi'(x_i)\geq \frac{x_i}{\log(x_i)}.$$
For each $i$ with $x_i$ large enough, \(x_i=2(s_i)_\R\) for some loxodromic element \(s_i \in \Gamma'\). Let us now consider 
$$m_i= \textrm{min}\{l \ |\  l \in \{2,5,8\}\ \textrm{satisfies}\  \textrm{Proposition~\ref{inedisp}} \mbox{ for } (s_i)_\R \}.$$
Furthermore, take $\tau_{m_i}$ and $\alpha_{m_i}$ as given in Lemma \ref{level_lemma} and Proposition \ref{inedisp}. Then the manifold $M_{i}=\Gamma_{\tau_{m_i}}(\alpha_{m_i})\backslash\mathbb{H}^n$ and the totally geodesic surface $S_i=\Gamma'_{\tau_{m_i}}(\alpha_{m_i})\backslash\mathbb{H}^2$ satisfy (see Remark \ref{rmk})
\begin{equation}\label{aprox}
\mathrm{Kiss}(M_i)\geq\mathrm{Kiss}(S_i) \geq \frac{x_i}{\log(x_i)}.
\end{equation}
On the other hand, take the isometry group $G_i=\Gamma/\Gamma_{\tau_{m_i}}(\alpha_{m_i})$ acting on the set of closed geodesics of $M_i$. By Lemma~\ref{C051121}(1), if we denote by $\gamma_1,\dots , \gamma_{\mathrm{Kiss}(S_i)}$ the systoles of $S_i$ embedded in $M_i$, the orbit sets $G_i \gamma_j$, $j\in\{1,\dots, \mathrm{Kiss}(S_i)\}$ are pairwise disjoint. It follows that

\begin{align}
\mathrm{Kiss}(M_i)&\geq \sum_{j=1}^{\mathrm{Kiss}(S_i)}|G_i\gamma_j|\\
&=\sum_{j=1}^{\mathrm{Kiss}(S_i)} \dfrac{|G_i|}{|(G_i)_{\gamma_j}|}\label{eq:aprox2},
\end{align}
where $(G_i)_{\gamma_j}$  denotes the isotropy group of $\gamma_j$ under the action of $G_i$. 
 By Lemma \ref{C051121}(2), $|(G_i)_{\gamma_j}|$ is at most the order of $\gamma_j$, and this is smaller than a fixed constant $C>0$ since $m_i\leq8$.
 Therefore, we get from \eqref{aprox} and \eqref{eq:aprox2} that
\begin{equation*}
\mathrm{Kiss}(M_i)\geq C \mathrm{Kiss}(S_i)\cdot|G_i| \geq C \frac{x_i}{\log(x_i)} |G_i|.
\end{equation*}
Since $\vol(M_i)=[\Gamma:\Gamma_{\tau_{m_i}}(\alpha_{m_i})]\vol(\Gamma\backslash\mathbb{H}^n)=|G_i|\vol(\Gamma\backslash\mathbb{H}^n)$, the above inequality becomes
\begin{equation}\label{aproxkiss}
\K(M_i)\geq \frac{C}{\vol(\Gamma\backslash\mathbb{H}^n)} \frac{x_i}{\log(x_i)} \vol(M_i).
\end{equation}

The goal now is to bound $\frac{x_i}{\log(x_i)}$ from below in terms of $\vol(M_i)$. Since $\Gamma_{\tau_{m_i}}(\alpha_{m_i})$ has index two in $\Gamma(\alpha_{m_i})$,  then
\begin{equation}\label{isotropy}
\vol(M_i)=[\Gamma:\Gamma_{\tau_{m_i}}(\alpha_{m_i})]\vol(\Gamma\backslash\mathbb{H}^n)=\frac{\vol(\Gamma\backslash\mathbb{H}^{n})}{2} [\Gamma:\Gamma(\alpha_{m_i})].
\end{equation}
Note that the norm of the ideal $(\alpha_{m_i})$ goes to infinity (see Lemma \ref{level_lemma}), thus we can apply the results in \cite[Section 5]{Mur19} along with the fact that $m_i\leq 8$, to obtain the following 

\begin{equation*}
\begin{aligned}
[\Gamma:\Gamma(\alpha_{m_i})]\leq \N(\alpha_{m_i})^{\frac{n(n+1)}{2}}&=O\left(|\alpha_{m_i}|^{\frac{n(n+1)}{2}}\right) \\
&=O\left(\left( x_i^{\frac{2(m_i+1)}{3}}\right)^{\frac{n(n+1)}{2}}\right)\\
&=O\left( x_i^{3n(n+1)}\right).
\end{aligned}
\end{equation*}

By putting together the estimate obtained above with Equation~\eqref{isotropy}, we find the lower bound 
\begin{equation}\label{volu}
x_i\geq C_1\cdot\vol(M_i)^{\frac{1}{3n(n+1)}}
\end{equation}
for some constant $C_1>0$. Finally, since the function $x\mapsto \frac{x}{\log(x)}$ is increasing for large $x$, by ~\eqref{aproxkiss} we get that 
\begin{equation*}
\K(M_i)\geq C_2\  \frac{\vol(M_i)^{1+\frac{1}{3n(n+1)}}}{\log(\vol(M_i))},
\end{equation*}
for a constant $C_2>0$. This ends the proof of Theorem \ref{maintheorem}.
\qedhere
\end{proof}

%% file: Sec7_dimension_three.tex
\subsection{Hyperbolic 3-manifolds}
In the three-dimensional case it will be convenient to consider the upper-half model of the hyperbolic $3$-space given by
 $$\mathbb{H}^3=\lbrace (z,t)\in\mathbb{C}\times\mathbb{R}; t>0\rbrace,$$ with the Riemannian metric $ds^2=\frac{dz^2+dt^2}{t^2}$. Note that we realize $\mathbb{H}^3$ as a subset of Hamilton's quaternion algebra
$$\mathcal{H}=\lbrace a+bi+cj+dk|a, b, c ,d\in\mathbb{R}, i^2=j^2=k^2=-1 \rbrace,$$ 
where we represent a point $P\in \mathbb{H}^3$ as a Hamiltonian quaternion $$P=(z,t):=x+yi+tj,$$ where $z=x+iy$.
Moreover $G=\SL(2,\mathbb{C})$ acts by isometries on $\mathbb{H}^3$, and this action is described as follows. For each $M=\left(\begin{array}{cc}
	a & b \\
	c & d \\
	\end{array}\right)\in \SL(2,\mathbb{C})$
	$$P\mapsto MP:=(aP+b)(cP+d)^{-1},$$
	where the inverse is taken in the skew field of Hamilton's quaternions.
This action is not faithful since \(-\mathrm{I}\) acts trivially, but the finite quotient $\mathrm{PSL(2,\mathbb{C})}=\SL(2,\mathbb{C})/\{\pm \mathrm{I}\}$ is isomorphic to $\mathrm{Isom}^{+}(\mathbb{H}^3)$ (see \cite[Chapter. 1]{EGM13}). As we have already observed with Spin groups, there is no loss of generality in identifying elements in \(\SL(2,\CC)\) with their projection in \(\PSLC\).
\label{typesofelements}
We recall that an element $\gamma\in \mathrm{SL}(2,\mathbb{C})$ is said to be:
	\begin{itemize}
	\item {\it elliptic} if $\gamma$ is conjugate to  $\left(\begin{array}{cc}
	\eta & 0 \\
	0 & \eta^{-1} \\
	\end{array}\right)$, $\mbox{ with }|\eta|=1, \eta \neq \pm 1.$\\	
	 \item {\it parabolic} if $\gamma$ is conjugate to  $\left(\begin{array}{cc}
		1 & z \\
		0 & 1 \\
		\end{array}\right)$, $\mbox{ for some }z\in\mathbb{C},\ z \neq 0.$\\
	 \item {\it loxodromic} if $\gamma$ is conjugate to  $\left(\begin{array}{cc}
	re^{i\theta} & 0 \\
	0 & r^{-1}e^{-i\theta} \\
	\end{array}\right)$, $r>1, r,\theta \in \mathbb{R}.$
				
	\end{itemize}
				
We define the trace of $\gamma=\pm \left(\begin{array}{cc}
	a & b \\
	c & d \\
	\end{array}\right)	\in \mathrm{PSL}(2,\mathbb{C})$
	as
	$$\tr(\gamma):=\mu \cdot (a+d),$$

where $\mu \in \{\pm 1\}$ is chosen so that $\tr(\gamma)=re^{i\theta}$ with $r\geq0$ and $\theta \in [0, \pi)$.

Consider the eigenvalues of $\gamma$ as the eigenvalues of a lift to $\mathrm{SL}(2,\mathbb{C})$. Hence the roots of the characteristic polynomial relative to $\gamma$ are  
	$$\lambda^{\pm }_{\gamma}=\frac{\tr(\gamma)\pm \sqrt{(\tr(\gamma))^2-4}}{2}.$$

These are the eigenvalues of one of the representatives of $\gamma \in \mathrm{SL}(2,\mathbb{C})$, and we denote by $\lambda_{\gamma}$ the eigenvalue with norm greater than one. We also choose the branch of the argument function $\Arg(z)$ on $\CC \setminus (-\infty,0]$ with $\Arg(z) \in (-\pi,\pi)$. It is well-known that $\lambda_{\gamma}$ determines the translation length $\ell(\gamma)$ of a  $\gamma$. More precisely, 
	\begin{equation}\label{leng}
	\ell(\gamma)=2\log(|\lambda_{\gamma}|) .    
	\end{equation}
					
	The \textit{holonomy} of $\gamma$ is defined as
	\begin{equation}\label{holo}
	 \theta(\gamma):=2 \mathrm{Arg}(\lambda_{\gamma}),   
	\end{equation}
					
We end this section by recalling how $l(\gamma)$ can be determined
from $\tr(\gamma)$. 
	\begin{prop}\label{traceequation}
	For any loxodromic element $\gamma \in \mathrm{SL}(2,\mathbb{C})$ we have
   $$4\cosh\left(\frac{l(\gamma)}{2}\right)=|\tr(\gamma)-2|+|\tr(\gamma)+2|.$$
   In particular $$4\cosh(l(\gamma))=|\tr(\gamma)^2|+|\tr(\gamma)^2-4|.$$
	\end{prop} 
	
\begin{proof}
See \cite[Proposition. 2.1]{DM21}, \textit{c.f}  \cite[Lemma. 5.1]{Gendulphe}. 
\end{proof}

\subsection{Arithmetic Kleinian groups}

A {\it Kleinian} group is a discrete group of $\PSLC$. Let $k$ be a number field with exactly one complex place and let $\mathcal{A}$ be a quaternion algebra over $k$ ramified at all real places. 
A Kleinian group $\Gamma$ is \emph{arithmetic} if it is commensurable with the projection  $P\rho(\mathcal{O}^1):=\rho(\mathcal{O}^1)/\{\pm \mathrm{I}\}$, where $\rho$ be a $k$-embedding of $\mathcal{A}$ into $M_2(\mathbb{C})$ and $\mathcal{O}^1$ denotes the group of elements of reduced norm one of an order $\mathcal{O}$ of $\mathcal{A}$. 
When $\Gamma\subset P\rho(\mathcal{O}^1)$ we say that $\Gamma$ is \emph{derived from a quaternion algebra}.

A hyperbolic $3$-orbifold $\Gamma\backslash\mathbb{H}^3$ is \emph{arithmetic} if $\Gamma$ is an arithmetic Kleinian group. For explicity examples, consider for each square-free positive integer $d>0,$ the Kleinian group $\PSL(2,\mathcal{O}_d)$, where $\mathcal{O}_d$ is the ring of integers of \(\Q(\sqrt{-d}).\) These groups are known as \emph{Bianchi groups} (see \cite[Chapter 7]{EGM13}).

Let $\mathcal{A}$ be a quaternion algebra over a number field $k$, with ring of integers $\mathcal{O}_k$. For any ideal $I \subset \mathcal{O}_k$, and any order $\mathcal{O} \subset \mathcal{A}$ we have an ideal $I\mathcal{O}$ defined by $$I\mathcal{O}=\left\{\sum_{j}t_jw_j| \ t_j\in I \ \textrm{and} \ w_j \in \mathcal{O} \right\}.$$ \indent The principal congruence subgroup of $\mathcal{O}^1$ of level $I$ is then given by
$$ \mathcal{O}^1(I)=\{\gamma\in \mathcal{O}^1 \ |\ \gamma-1 \in I\mathcal{O}\}.$$

Suppose that $\Gamma < \PSLC$ is a Kleinian group. We will denote by $\tilde{\Gamma}$ the preimage of $\Gamma$ by the natural projection of $\SL(2,\mathbb{C})$ into $\PSL(2,\mathbb{C})$. In particular, if $\Gamma$ is arithmetic derived from a quaternion algebra, then $\tilde{\Gamma} < \mathcal{O}^1$ for some order $\mathcal{O}$. In this way, for any ideal $I$ we define the principal congruence subgroup of $\Gamma$ of level $I$ as the projection of $\tilde{\Gamma} \cap \mathcal{O}^1(I)$ onto $\Gamma$. If $I$ is a principal ideal, generated by $\alpha\in\mathcal{O}_k$, we denote $\Gamma(I)$ by $\Gamma(\alpha)$ instead of $\Gamma(\langle\alpha\rangle)$. The key fact about congruence subgroups is the following lemma (Compare\ with\  Lemma~\eqref{realpart}). 
\begin{lemma}\label{tracelevel}
For any $\gamma \in \mathcal{O}^1(I)$ we have $\tr(\gamma)\equiv 2\ \mathrm{mod}\ I^2.$
\end{lemma}
\begin{proof}
Let \(\gamma \in \mathcal{O}^1(I)\), by definition we can write \(\gamma =1+\eta\), with \(\eta \in I\mathcal{O}\). Since \(\tr(\eta)=\eta+\eta^* \in I \) and \(\mathrm{N_{red}}(\gamma) =\eta \eta^*\in I^2\) (see \cite[Lemma 3.3]{KS}), we have
\[1=\gamma\gamma^*=1+\tr(\eta)+ \mathrm{N_{red}}(\gamma). \]
Therefore \(\tr(\eta) \in I^2\), and then \(\tr(\gamma) \equiv 2\  \mathrm{mod}\ I^2\). 
\end{proof}
\subsection{Displacement estimates for congruence subgroups}
We will now construct hyperbolic 3-orbifolds for which we can determine their set of systoles. Before that, we introduce some notation that will be convenient for this purpose.

Let $T:\CC^* \to \CC^*$ be the surjective holomorphic map given by $T(z)=z+z\inv$, and let $V=\CC \setminus (-\infty,0]$. 
We observe that, if $\gamma\in\PSLC$ is loxodromic, and $t=\tr(\gamma)$ with largest eigenvalue $\lambda$, then $T(\lambda)=t$.
On the other hand, $\Arg(T(z))$ is a continuous map from $T\inv(V)$ to $(-\pi,\pi)$. If we write $z$ in its polar coordinates 
with $|\Arg(z)|<\pi$, we have that \[T(z)=(|z|+|z|\inv)\cos(\Arg(z))+i(|z|-|z|\inv)\sin(\Arg(z)).\]
Thus,  $|\Arg(T(z))|=\frac{\pi}{2}$ if, and only if, $|\Arg(z)|=\frac{\pi}{2}$. Moreover, if $|z|>1$ and $|\Arg(z)| \neq \frac{\pi}{2}$ we obtain 
\begin{equation}\label{lemma:trace-eingenvalue}
\tan(\Arg(T(z))) = \frac{|z|-|z|\inv}{|z|+|z|\inv} \tan(\Arg(z)).  
\end{equation}

\begin{prop}\label{prop:elements_square_systoles}
 Let $k$ be an imaginary quadratic field, and $\Gamma$ an arithmetic Kleinian group derived from a quaternion algebra over $k$. Then, there exist $L,\epsilon>0$ such that if $ \gamma \in \Gamma$ is a loxodromic element with $$\ell(\gamma) > L \mbox{ and } 0 \leq \hol(\gamma) < \epsilon,$$ then $\gamma^2$ realizes the systole of $\Gamma(\tr(\gamma)) \backslash \Hyp$.
\end{prop}

\begin{proof}
By hypothesis, there exist a quaternion algebra $A$ over $k$, an order $\mathcal{O}$ of $A$ and a $k$-monomorphism of algebras $\rho: A \to \SL(2,\CC)$ such that $\tilde{\Gamma} < \rho(\mathcal{O}^1)$. Hence, we can assume that $\tilde{\Gamma} < \mathcal{O}^1.$

Suppose that $\tr(\gamma)=t\in \mathcal{O}_k$, and let $\tilde{\gamma}$ be a representative of $\gamma$ in $\tilde{\Gamma}$. By definition  $$\tilde{\Gamma}(t)=\langle \tilde{\Gamma} \cap \mathcal{O}^1(t), -1 \rangle,$$ and $\Gamma(t)=\tilde{\Gamma}(t)/\{\pm 1\}$. Since $\tilde{\gamma}^2-t\tilde{\gamma}+1=0$ in $\mathcal{O}$ we get that $-\tilde{ \gamma}^2 \in \tilde{\Gamma}(t)$, and then $\gamma^2 \in \Gamma(t)$. 

Now, let $\lambda\in \CC$ be the eigenvalue of $\gamma$ with $|\lambda|>1$. By ~\eqref{leng} we have $\ell(\gamma^2)=2\ell(\gamma)=4\log(|\lambda|)$. Since $\gamma^2\in\Gamma(t)$ we get that
 $$\sys(\Gamma(\tr(\gamma)) \backslash \Hyp) \leq 4 \log(|\lambda|)=\ell(\gamma^2).$$
 
Our purpose is to give conditions on \(t\) such that this inequality becomes an equality.
Let $\eta \in \Gamma(t)$ be a loxodromic element. There exists a representative, say $\tilde{\eta}$, of $\eta$ with $\tilde{\eta} \in \tilde{\Gamma} \cap \mathcal{O}^1(t)$. It follows from Lemma \ref{tracelevel} that, if $\tau=\tr(\eta)$, then $\tau=2+t^2\zeta$ for some $\zeta \in \mathcal{O}_k$ and $\zeta \neq 0$. If $\tau=T(\mu)$ with $|\mu|>1$, since $t^2=\lambda^2+\lambda^{-2}+2$, then $\tau=2+t^2\zeta$ can be rewritten as
 \begin{equation}\label{eq:relation_T}
 T(\mu)=2(\zeta+1)+\zeta T(\lambda^2).
 \end{equation} 
 
By \eqref{leng} it is sufficient to show that  $|\mu|>|\lambda|^2$. We will divide our analysis in two cases:
 
\textit{Case 1.} $\zeta \notin \mathcal{O}^*_k$: Since $k$ is a quadratic field, we have $|\zeta|^2 \geq 2$, and then $|\zeta| \geq \sqrt{2}$. Firstly, we can rewrite \eqref{eq:relation_T} as 
\begin{align}\label{eq:relation_T2}
	T(\mu)&=\zeta\lambda^2\left(1+(\lambda^{2})^{-2}+2(\lambda^{2})\inv+2(\lambda^{-2} \zeta\inv)\right)\\
	&=\zeta\lambda^2R(\lambda^2,\zeta),
\end{align}
where $R(z,\theta)=1+z^{-2}+2z\inv+2z\inv \theta\inv$ is defined on $\CC^* \times (\mathcal{O}-\{0\})$. Since $|\theta|\geq 1$ for any $\theta \in \mathcal{O}-\{0\}$, it follows that for any $\delta >0$ there exists $N>0$ such that for $(z,\theta)\in\CC^* \times (\mathcal{O}-\{0\})$ with $|z|>N$, it holds $|R(z,\theta)|>1-\delta.$ In particular, we can choose \(N_0>2\) such that $|z|>N_0$ implies $|R(z,\theta)|>\dfrac{3\sqrt{2}}{4}.$

\textit{Case 2.} $\zeta \in \mathcal{O}^*_k$: Since \(k\) is a quadratic imaginary field, \(\zeta\) is a \(n \mathrm{th}\) root of unity. In particular, \(\zeta\) is an algebraic integer of degree \(\phi(n)\), where \(\phi\) denotes the Euler's totient function (see \cite[Chapter IV, Theorem 2]{Lang13}). Since \(\zeta\) has degree \(2\),  we conclude that \(n \in \{1,2,3,4,6\} \), i.e.  $\zeta \in J=\{\pm 1, \pm i, \pm \omega, \pm \omega^2\}$, where $\omega=\frac{1}{2}-i\frac{\sqrt{3}}{2}$ is a primitive sixth root of unity. 

By Proposition \ref{traceequation} we have that
\[4\cosh\left(\frac{\ell(\eta)}{2}\right)=|T(\mu)-2|+|T(\mu)+2|.\]

However, $|\zeta|=1$ and then Equation \eqref{eq:relation_T} implies that
\[|T(\mu)-2|+|T(\mu)+2|=|T(\lambda^2)+2|+|T(\lambda^2)+2+4\zeta\inv|.\]

Hence, we can guarantee that $\ell(\eta) \geq 2\ell(\gamma)$ whenever 
\begin{equation}\label{eq:cond_lambda}
|T(\lambda^2)+2+4\zeta\inv| \geq |T(\lambda^2)-2|
\end{equation}

for any $\zeta\in J$. In order to compute the difference $|T(\lambda^2)+2+4\zeta\inv|-|T(\lambda^2)-2|$, we consider the map defined on $(-\pi/2,\pi/2)$, by
$$h_{P,\zeta}(\phi)=|Pe^{i\phi}+2+4\zeta|^2-|Pe^{i\phi}-2|^2$$
where $P>1$ and $\zeta\in\CC$ are fixed. In this way Inequality \eqref{eq:cond_lambda} is equivalent to
\begin{equation}\label{eq:h_P}
h_{|T(\lambda^2)|,\zeta^{-1}}(\Arg(T(\lambda^2)))\geq 0
\end{equation}  
for any $\zeta\in J$. We then look for conditions on $\Arg(\lambda)$ such that \eqref{eq:h_P} holds for any $\zeta\in J$. It is clear that $h_{P,-1}(\phi)=0$ for any $P$ and $\phi$, so we can assume that $\zeta\in J\setminus\{-1\}$.

It is straightforward to check that 
\[h_{P,\zeta}(\phi)=16(1+\Re(\zeta))+8P\cos(\phi)[1+\Re(\zeta)+\Im(\zeta)\tan(\phi)]\].

Hence, if \(\zeta=1\), since \(\cos(\phi)>0\) we have \(h_{P,1}(\phi)>0\) for all \(\phi \in (-\pi/2,\pi/2) \). It follows from $\zeta\in J\setminus\{\pm 1\}$ that
\begin{equation}\label{eq:part_real_zeta}
1+\Re(\zeta) \geq \frac{1}{2}
\end{equation}
and
\begin{equation}\label{eq:part_im_zeta}
\Im(\zeta) \geq -1 \mbox{ , } \Im(\zeta) \neq 0.
\end{equation}

Therefore $h_{P,\zeta}>0$ whenever 
\begin{equation}\label{eq:positivity}
1+\Re(\zeta)+\Im(\zeta)\tan\phi>0.
\end{equation} 

Suppose now that $0<\tan(\phi)<\frac{1}{2}$. If $\Im(\zeta)\geq 0$ then \eqref{eq:positivity} follows from \eqref{eq:part_real_zeta}. On the other hand, if $\Im(\zeta)<0$, by \eqref{eq:part_im_zeta} we get $0<-\Im(\zeta)\leq1$, and together with \eqref{eq:part_real_zeta} we obtain that 

\[\tan(\phi)<\frac{1}{2}<\frac{1+\Re(\zeta)}{-\Im(\zeta)},\]
from which \eqref{eq:positivity} follows. Then, if $0<\Arg(\lambda)<\frac{1}{2}\arctan\left(\frac{1}{2}\right)$, then $0<\tan(\Arg(T(\lambda^2)))<\frac{1}{2}$ by \eqref{lemma:trace-eingenvalue} and the fact that $\Arg(\lambda^2)=2\Arg(\lambda)$. Therefore  \eqref{eq:h_P} follows as desired.

Since $\ell(\gamma^2)=2\ell(\gamma)$ and $\hol(\gamma^2)=2\hol(\gamma)=2\Arg(\lambda)$, we conclude from the analysis of the two cases for $L=4\log(N_0)>0$ (with $N_0$ given in Case 1), and $\epsilon=\frac{1}{2}\arctan\frac{1}{2}$, that if $\ell(\gamma)>L$ and $0\leq\hol(\gamma)<\epsilon$, then $\gamma^2$ minimizes the set of displacements of $\Gamma(\tr(\gamma))$, and therefore
\[\sys(\Gamma(t) \backslash \Hyp) = 2\ell(\gamma).\]
\end{proof}

\subsection{Proof of Theorem \ref{secondresult}}

Let $\gamma \in \PSL(2,\CC)$ be a loxodromic element. We can associate to $\gamma$ the complex number $z(\gamma)=e^{\frac{\ell(\gamma)}{2}}e^{i\hol(\gamma)}$. 
Thus, by construction in Section \eqref{typesofelements} we have that $T(z(\gamma))$ is the trace of some lifting of $\gamma$ in $\SL(2,\CC)$.
Thereby, we will adopt as the \emph{trace} of $\gamma$ the complex number written as $$T(\gamma):=T(z(\gamma)).$$

Note that this definition of trace remains invariant under conjugation and  extend the definition of trace to the conjugacy class of any subgroup of $\PSL(2,\CC)$. For a complex number $z$ we will define the \emph{norm of $z$} as the nonnegative real number $|z|^2$.  
If $\Gamma<\PSL(2,\CC)$ is a Kleinian group, we define $\sigma(N,I)$ (resp. $\tau(N,I)$) as the number of pri\-mitive conjugacy classes of $\Gamma$ with norm of trace at most $N$ and holonomy in $I,$ counted with multiplicity (resp. counted without multipliticy). By definition, the {\it mean multiplicity} is given by 
 $$\mu_0(N,I)=\dfrac{\sigma(N,I)}{\tau(N,I)}.$$ 

These definitions will be convenient for presenting the following proposition.

\begin{prop}\label{prop:mean_multiplicity_traces}
 Let $\Gamma$ be an arithmetic Kleinian group derived from a quaternion algebra over an imaginary quadratic field $k$. For any subinterval $I \subset [0,2\pi]$, let $\mu_0(N,I)$ be the mean multiplicity of primitive conjugacy classes of $\Gamma$ with trace of norm at most $N$ and holonomy contained in $I$. Then there exists a constant $c>0$ depending only on $k$ and $I$ such that
 \[ \mu_0(N,I) \gtrsim c\frac{N}{\log(N)} \mbox{ when } N \to \infty .\] 
 
\end{prop}
\begin{proof}

Let $\ok$ be the ring of integers of $k$. For any conjugacy class $[\gamma] \subset \Gamma$  we have $T(\gamma) \in \ok$. Moreover, if $|z(\gamma)|>1$, then the norm of $T(\gamma)$ is at most $|z(\gamma)|^2+3$.
For any $L>0$ and $I \subset [0,2\pi]$, consider 
\[\nm(L,I)=\# \{[\gamma] \subset \Gamma \mid \gamma \mbox{ is primitive, } |z(\gamma)| \leq L \mbox{ and } \hol(\gamma) \in I\}. \] 
By \cite[Corollary]{SW99} (see also \cite[Thm. 1.3]{Marg14} for a more explicit statement), there exists a constant $c_1$ which depends only on $I$ such that
\[\nm(L,I) \sim c_1\frac{L^4}{\log(L)} \mbox{ when } L \to \infty. \]

Hence,
$\sigma(N,I)$ is at least $\nm(\sqrt{N-3},I)$, implying that
\begin{equation}\label{eq:Margulis_grow}
\sigma(N,I) \gtrsim c_1' \frac{N^2}{\log(N)} 
\end{equation}
for some constant $c_1$ depending only on $I$ and $N$ sufficiently large.

On the other hand, since $\ok$ is a lattice in $\CC$, there exists a constant $c_2>0$ depending only on $k$ such that 
\[\#(\ok \cap B(0,R)) \sim c_2R^2 \]
(see \cite[Ch. V, Thm. 2]{Lang13}). Hence, when $N$ is big enough we have
\begin{equation}\label{eq:lattice_grow}
\tau(N,I) \lesssim c_2N 
\end{equation}

Thus, if we combine the asymptotic bounds \eqref{eq:Margulis_grow} and \eqref{eq:lattice_grow}, by the definition of mean multiplicity
there exists a constant $c>0$ which depends only in $I$ and $k$ such that
\[\mu_0(N,I) \gtrsim c \frac{N}{\log(N)} \text{ when } N \to \infty.\]

\end{proof}

We are now ready to finish the proof of Theorem \ref{secondresult}. In fact, we will state a more precise result which implies that every commensurability class of arithmetic hyperbolic $3$-manifolds with imaginary quadratic invariant trace field contains a sequence of manifolds with kissing number as stated. It is well known by the Classification Theorem of Quaternions Algebras over number fields (see \cite[Theorem 7.3.6]{MR}) that there exist compact and non compact arithmetic hyperbolic $3-$manifolds with this property.
\begin{theorem}
 Let $\Gamma < \PSLC$ be an arithmetic Kleinian group with 
 an imaginary  quadratic invariant trace field. There exists a sequence $\{\Gamma_j\}$ of torsion-free subgroups of $\Gamma$ with arbitrarily large index such that the corresponding sequence of finite volume hyperbolic 3-manifolds $M_j=\Gamma_j \backslash \Hyp$ satisfies
 \[\K(M_j) \geq c \frac{\vol(M_j)^\frac{4}{3}}{\log(\vol(M_j)} \]
 for some constant $c$ which does not depend on $j$.
\end{theorem}
\begin{proof}
We can suppose that $\Gamma$ is derived from a quaternion algebra since $\Gamma^{(2)}=\langle \gamma^2 \mid \gamma \in \Gamma \rangle$ has finite index in $\Gamma$ and is derived from a quaternion algebra (\cite[Corollary 8.3.5]{MR}).

Consider the constants $L$ and $\varepsilon$ given in Proposition \ref{prop:elements_square_systoles}, and let $k$ be the invariant trace field of $\Gamma$. If we set $I=[0,\varepsilon]$, then by Proposition \ref{prop:mean_multiplicity_traces} there exists a sequence of traces $t_j \in \ok$ with $|t_j|^2=N_j \to \infty$ such that the number $n_j$ of primitive conjugacy classes of $\Gamma$ with trace $t_j$ satisfy \[n_j \geq c \dfrac{N_j}{\log(N_j)},\]  where $c=c_{k,\varepsilon}>0$. 

By Lemma \eqref{tracelevel}, if \(N_j\) is large enough, \(\tr(\gamma)=\pm 2\) or \(|\tr(\gamma)|>2\). Thus we can assume that for \(j\) sufficiently large $\Gamma(t_j)$ is torsion-free and $N_j>L$ for all $j$. Now we can argue as in theorem \ref{maintheorem}. 

Let $\gamma_1,\ldots,\gamma_{n_j} \in \Gamma$ be a set of representatives of all primitive conjugacy classes in $\Gamma$ of trace $t_j$, and let $G_j=\Gamma / \Gamma(t_j)$ be a subgroup of isometries of $M_j=\Gamma(t_j) \backslash \Hyp$. By Proposition \ref{prop:elements_square_systoles}, each $\gamma_i^2 \in \Gamma(t_j)$, and their induced closed geodesics are systoles of $M_j$.  Since $\gamma_i^2$ has order $2$, by Lemma \ref{C051121} we have that 
\[\K(M_j) \geq \sum_{i=1}^{n_j} \# (G_j \cdot \gamma_i^2) \geq \frac{1}{2} n_j\# G_j  .\]

It is a well-known fact that $\# G_j = [\Gamma: \Gamma(t_j)] \leq CN_j^3$ for some constant $C>0$, which does not depend on $j$ (see \cite{Mak},\cite{KS} or \cite{Kuc}). Moreover, $\vol(M_j)=\mu [\Gamma:\Gamma(t_j)]$ for all $j$, where $\mu=\vol(\Gamma \backslash \Hyp)$ . Putting together the inequalities above, as in the proof of Theorem \ref{maintheorem}, we obtain that
\[ \K(M_j) \geq c \frac{\vol(M_j)^{\frac{4}{3}}}{\log(\vol(M_j))} \]
for any $j$, where $c$ does not depend on $j$.
 
\end{proof}